# Recursive estimation of possibly misspecified MA(1) models: Convergence of a general algorithm

James L. Cantor[1] and David F. Findley[2]

*Science Application International Corporation and U.S. Census Bureau*

**Abstract:** We introduce a recursive algorithm of conveniently general form for estimating the coefficient of a moving average model of order one and obtain convergence results for both correct and misspecified MA(1) models. The algorithm encompasses Pseudolinear Regression (PLR—also referred to as AML and $RML_1$) and Recursive Maximum Likelihood ($RML_2$) without monitoring. Stimulated by the approach of Hannan (1980), our convergence results are obtained indirectly by showing that the recursive sequence can be approximated by a sequence satisfying a recursion of simpler (Robbins-Monro) form for which convergence results applicable to our situation have recently been obtained.

## 1. Introduction and overview

Our focus is on estimating the coefficient $\theta$ of an invertible scalar moving average model of order 1 (MA(1)),

$$y_t = \theta e_{t-1} + e_t \tag{1.1}$$

where $e_t$ is treated as an unobserved, constant-variance martingale-difference process. We do not assume the series $y_t, -\infty < t < \infty$ from which the observations come is correctly modeled by (1.1). They can come from any invertible autoregressive moving average (ARMA) model or from more general models; see Section 2. What we seek is a $\theta$ that minimizes the loss function

$$\bar{L}(\theta) = E[(y_t - y_{t|t-1}(\theta))^2] = E[e_t^2(\theta)] \tag{1.2}$$

where $e_t(\theta) = y_t - y_{t|t-1}(\theta)$ and $y_{t|t-1}(\theta)$ is the one-step-ahead-prediction of $y_t$ from $y_s, -\infty < s \leq t-1$ based on the model defined by $\theta$ (see (2.7) below). We define *optimal* estimation procedures to be those whose sequence of estimates $\theta_t$ minimizes (1.2) in the limit. This is a property of (nonrecursive) maximum likelihood-type estimates of $\theta$, see Pötscher [23].

In this article, we analyze a continuously indexed family of recursive procedures for estimating $\theta$. Recursive procedures form an estimate $\theta_t$ for time $t$ using the observation $y_t$ at time $t$, the estimate $\theta_{t-1}$ for $t-1$ and other recursively defined quantities. Our family encompasses two standard algorithms, Recursive Maximum

---

[1]Science Applications International Corporation (SAIC), 4001 North Fairfax Drive, Suite 250, Arlington, VA 22203, e-mail: `james.l.cantor@saic.com`

[2]U.S. Census Bureau, Statistical Research Division, Room 3000-4, Washington, DC 20233-9100, e-mail: `david.f.findley@census.gov`

*AMS 2000 subject classifications:* primary 62M10; secondary 62L20.

*Keywords and phrases:* time series, Robbins-Monro, PLR, AML, $RML_1$, $RML_2$, misspecified models.

Likelihood (RML) which is referred to throughout as RML$_2$ [12, 21], and the simpler Pseudolinear Regression (PLR) [21]—also known as Approximate Maximum Likelihood (AML) [24] and RML$_1$ [11, 20]. More specifically, our general recursive algorithm generating $\theta_t$ depends on an index $\beta, 0 \leq \beta \leq 1$. The algorithm reduces to PLR when $\beta = 0$ and to RML$_2$ when $\beta = 1$.

Our main convergence result, Theorem 4.1, is obtained by constructing an approximating sequence $\hat{\theta}_t$ for which $\theta_t - \hat{\theta}_t \xrightarrow{a.s.} 0$ holds and which satisfies a Robbins-Monro recursion,

$$(1.3) \qquad \hat{\theta}_t = \hat{\theta}_{t-1} - \delta_t f(\hat{\theta}_{t-1}, \beta) + \delta_t \gamma_t,$$

in which $\gamma_t \xrightarrow{a.s.} 0$ and $\delta_t > 0, \delta_t \xrightarrow{a.s.} 0, \sum_{k=0}^{\infty} \delta_k = \infty$ a.s., and

$$(1.4) \qquad f(\theta, \beta) = -\int_{-\pi}^{\pi} \frac{e^{i\omega} + \beta\theta}{|(1 + \theta e^{i\omega})(1 + \theta\beta e^{i\omega})|^2} g_y(\omega) d\omega.$$

Here $\xrightarrow{a.s.}$ denotes almost sure convergence (convergence with probability one) and $g_y(\omega)$ denotes the spectral density of the time series $y_t$. Note that when $\beta = 0$, then

$$(1.5) \qquad f(\theta, 0) = -\int_{-\pi}^{\pi} \frac{e^{i\omega}}{|(1 + \theta e^{i\omega})|^2} g_y(\omega) d\omega = -E[e_{t-1}(\theta) e_t(\theta)],$$

and when $\beta = 1$, then

$$(1.6) \qquad f(\theta, 1) = -\int_{-\pi}^{\pi} \frac{e^{i\omega} + \theta}{|(1 + \theta e^{i\omega})^2|^2} g_y(\omega) d\omega = \frac{1}{2} \frac{d}{d\theta} E[e_t^2(\theta)] = \frac{1}{2} \bar{L}'(\theta)$$

where $\bar{L}'(\theta)$ denotes the first derivative of $\bar{L}(\theta)$. We then apply a result of Fradkov implicit in [8], as extended and corrected by Findley [9], to show that $\hat{\theta}_t$ converges to $\{\theta \in \Theta \colon f(\theta, \beta) = 0\}$ where $\Theta$ is the open interval (-1,1) of real $\theta$ with $|\theta| < 1$. (A similar result is implicit in proofs of Theorems 2.2.2–2.2.3 of Chen [7].) Hence, for $\beta = 0$, $\theta_t \xrightarrow{a.s.} \{\theta \in \Theta \colon E[e_{t-1}(\theta) e_t(\theta)] = 0\}$ and for $\beta = 1$, $\theta_t \xrightarrow{a.s.} \{\theta \in \Theta \colon \bar{L}'(\theta) = 0\}$. Here and below, $\theta_t$ convergence a.s. to a set means that except on a set of $\xi \in \Xi$ with probability zero, every cluster point of $\theta_t(\xi)$ is an element of the set.

In the incorrect model situation, in which $g_y(\omega)$ is not proportional to $|1 + \theta e^{i\omega}|^2$, for examples we have analyzed [5], these zero sets will be disjoint, establishing that PLR converges to different values than RML$_2$. Consequently, under the assumptions of Theorem 4.1, we recover the results of Cantor [4] that were given in separate theorems and proofs, establishing that, for certain families of AR(1) and MA(2) processes, RML$_2$ estimates of $\theta$ in the model (1.1) converge to an optimal limit (a minimizer of (1.2)) whereas PLR estimates converge to a suboptimal limit [4, 5].

When the data come from an invertible MA(1) model, it is known that PLR and monitored versions of RML$_2$ can provide strongly consistent estimates of $\theta$ [4, 11, 17, 19]. More generally, in the correct model situation for ARMAX models, i.e., ARMA models with an exogenous input, Lai and Ying [17] provided a rigorous proof of strong consistency of PLR (under a positive real condition on the MA polynomial) and also of a monitored version of RML$_2$ whose monitoring scheme involves non-linear projections and an intermittently used recursive estimator for which consistency has already been established. In Section 4 of [19], Lai and Ying consider a simpler modification of RML$_2$ in which, for monitoring, only auxiliary consistent recursive estimates are used. They present detailed outlines of proofs of strong consistency and asymptotic normality of the estimates from this new monitored RML$_2$ scheme. The construction of Section II of [18] can be used to obtain auxiliary recursive estimates with the properties required.

There is a rather comprehensive theory of recursive estimation of autoregressive (AR) models, encompassing certain incorrect model situations for algorithms like PLR (see e.g., [6]). There are, however, no published convergence results with rigorous proofs for MA models in the incorrect model situation. Ljung's seminal work on the convergence of recursive algorithms [20, 21] mentions the incorrect model situation but provides only suggestive results (further discussed in Section 5).

This article has five sections. In Section 2, the assumptions on the data and some consequences for the MA(1) model are given. In Section 3, the general recursive algorithm is presented. The Convergence Theorem is stated and proved in Section 4. Required preliminary technical results are given in Section 4.1 and the proof of the theorem is provided in Section 4.2. Finally, Section 5 concludes the article with a brief discussion.

## 2. Assumptions

The observations $y_t, t \geq 1$ are assumed to come from a mean zero, covariance stationary scalar series, $y_t, -\infty < t < \infty$ defined on the probability space $(\Xi, \mathcal{F}, P)$. We use the following additional assumptions on the process $y_t$:

(D1) $y_1$ is nonzero with probability one; i.e., $P\{y_1^2 > 0\} = 1$.
(D2) The series has a linear representation

$$(2.1) \qquad y_t = \sum_{s=0}^{\infty} \kappa_s \epsilon_{t-s} \text{ such that } \kappa_0 = 1 \text{ and } \sum_{s=0}^{\infty} |\kappa_s| < \infty$$

in which $\kappa(z) = \sum_{s=0}^{\infty} \kappa_s z^s$ is nonzero for $|z| \leq 1$ and $\{\epsilon_t\}$ is a martingale-difference sequence (m.d.s.) with respect to the sequence of sigma fields $\mathcal{F}_t = \sigma(y_s, -\infty < s \leq t)$. Thus $E[\epsilon_t | \mathcal{F}_{t-1}] = 0$. By a result of Wiener [25, Theorem VI 5.2], $\kappa(z)^{-1} = \sum_{s=0}^{\infty} \beta_s z^s$ with $\sum_{s=0}^{\infty} |\beta_s| < \infty$, whence

$$(2.2) \qquad \epsilon_t = \sum_{s=0}^{\infty} \beta_s y_{t-s} \ (\beta_0 = 1).$$

(D3) The conditional variance $E[\epsilon_t^2 | \mathcal{F}_{t-1}]$ is constant almost surely; i.e., $E[\epsilon_t^2 | \mathcal{F}_{t-1}] = \sigma_\epsilon^2$ a.s. Equivalently, $E[\epsilon_t^2] = \sigma_\epsilon^2$ and $\epsilon_t^2 - \sigma_\epsilon^2$ is a m.d.s. with respect to the $\mathcal{F}_t$.
(D4) $\{\epsilon_t\}$ is bounded a.s.; $\sup_t |\epsilon_t| \leq K$ a.s. for some $K < \infty$.

From (D2)–(D3), the spectral density $g_y(\omega)$ can be expressed as

$$(2.3) \qquad g_y(\omega) = \frac{\sigma_\epsilon^2}{2\pi} \left|\kappa(e^{i\omega})\right|^2 \text{ where } \kappa(e^{i\omega}) = \sum_{j=0}^{\infty} \kappa_j e^{ij\omega},$$

and

$$(2.4) \qquad 0 < m \leq g_y(\omega) \leq M < \infty \text{ for all } -\pi \leq \omega \leq \pi$$

for positive constants $m$ and $M$. The series $y_t$ is an invertible ARMA process if and only if $\kappa(z)$ is a rational function.

Assumption (D4) is used extensively in the proof of the convergence theorem, Theorem 4.1, in Section 4.

Under (D2)–(D4), we can apply, for example, the First Moment Bound Theorem of Findley and Wei [10] to show that $t^{-1} \sum_{s=j+1}^{t}(y_s y_{s-j} - \gamma_j^y) \xrightarrow{a.s.} 0$. Hence, from the particular case $y_t = \epsilon_t$ in (2.1) and $j = 0$,

$$(2.5) \qquad t^{-1} \sum_{s=1}^{t} \epsilon_s^2 \xrightarrow{a.s.} \sigma_\epsilon^2.$$

We consider models for $y_t$ of the invertible, stationary first-order moving-average type (MA(1)) given by

$$(2.6) \qquad y_t = \theta e_{t-1} + e_t, \quad -\infty < t < \infty.$$

For a given coefficient $\theta$ such that $|\theta| < 1$, the difference equation (2.6) is satisfied with $e_t = e_t(\theta)$ given by the mean zero, covariance stationary one-step-ahead-prediction-error series,

$$(2.7) \qquad e_t(\theta) = (1 + \theta B)^{-1} y_t = \sum_{j=0}^{\infty} (-\theta)^j y_{t-j} = y_t - y_{t|t-1}(\theta),$$

from the MA(1) predictor $y_{t|t-1}(\theta) = -\sum_{j=1}^{\infty}(-\theta)^j y_{t-j}$, see (5.1.21) of [3]. Here $B$ is the backshift operator; i.e., $By_t = y_{t-1}$. The coefficient $\theta$ is referred to as the MA *coefficient*. Thus,

$$(2.8) \qquad y_t = e_t(\theta) + \theta e_{t-1}(\theta).$$

The infinite series in (2.7) converges in mean square and, from (D4) and the representation (2.1), also almost surely. Thus, $e_t(\theta)$ represents the optimal one-step-ahead-prediction-error process from the perspective of the model (2.6). The model (2.6) is correct if $e_t(\theta)$ coincides (a.s.) with the m.d.s. $\epsilon_t$ in (2.2), in which case $\beta_s = (-\theta)^s, k \geq 0$. Whether or not the model is correct for any $\theta$, forecast errors $e_t(\theta)$ appearing in loss functions such as (1.2) and elsewhere are calculated as in (2.7). We emphasize that (2.1) allows data processes far more general than MA(1) processes. In particular, the $z$-transform, $\sum_{s=0}^{\infty} \kappa_s z^s$ is not required to be rational. For example, time series conforming to the exponential models of Bloomfield [2] have non-rational $\kappa(z)$ without zeroes in $|z| \leq 1$.

Let $\Theta = (-1, 1)$. From (2.7), the spectral density of $e_t(\theta)$ is $g_e(\theta, \omega) = g_y(\omega) \cdot |1 + \theta e^{i\omega}|^{-2}$, so for $\bar{L}(\theta)$ defined by (1.2), we have

$$(2.9) \qquad \bar{L}(\theta) = \int_{-\pi}^{\pi} \frac{g_y(\omega)}{|1 + \theta e^{i\omega}|^2} d\omega.$$

By (2.4) and the continuity of $g_y(\omega)$, $\bar{L}(\theta)$ is positive, infinitely differentiable, and nonconstant on the interior of $[-1, 1]$, i.e., on $\Theta$, and infinite at the endpoints. Therefore it has a minimum value over $[-1, 1]$ and

$$(2.10) \qquad \Theta^* \equiv \left\{ \theta \in [-1, 1] \colon \theta = \arg\min_{\theta \in [-1,1]} \bar{L}(\theta) \right\},$$

is a subset of $[-K, K]$ for some $0 < K < 1$. Also $\Theta^* \subseteq \Theta_0^* = \{\theta \in \Theta \colon \bar{L}'(\theta) = 0\}$. We are interested in a.s bounded random recursive sequences $\theta_t = \theta_t(\xi)$ that converge

a.s. to $\Theta^*$ or at least to $\Theta_0^*$. If $\Theta_0^*$ contains only one point, $\theta_0^*$, then $\theta_t$ converges to $\theta_0^*$ a.s. Our results will establish convergence of the sequence of estimates $\theta_t$ defined by the general algorithm presented below to the set of zeroes of $f(\theta, \beta)$ defined by (1.4).

## 3. The general recursive algorithm

For $0 \leq \beta \leq 1$, we define a general recursion for estimating the MA coefficient $\theta$ of (1.1):

$$(3.1a) \qquad \theta_t = \theta_{t-1} + \bar{P}_t^{-1} \frac{1}{t} \phi_{t-1} e_t; \quad \theta_1 = 0, \, t \geq 2,$$

$$(3.1b) \qquad \bar{P}_t = \frac{1}{t} \sum_{s=1}^{t-1} \phi_s^2 = \bar{P}_{t-1} + \frac{1}{t}[\phi_{t-1}^2 - \bar{P}_{t-1}]; \quad \bar{P}_1 = 0; \, t \geq 2,$$

$$(3.1c) \qquad e_t = y_t - \theta_{t-1} e_{t-1}; \quad e_1 = y_1, \, t \geq 2,$$

$$(3.1d) \qquad \phi_t = x_t - \theta_{t-1} \phi_{t-1}; \quad \phi_1 = x_1, \, t \geq 2,$$

$$(3.1e) \qquad x_t = y_t - \beta \theta_{t-1} x_{t-1}; \quad x_1 = y_1, \, t \geq 2.$$

From (3.1a), it follows for $0 \leq s \leq t-1, t \geq 2$ that

$$(3.2) \qquad \theta_{t-s} = \theta_t - \sum_{l=0}^{s-1} (t-l)^{-1} \bar{P}_{t-l}^{-1} \phi_{t-l-1} e_{t-l},$$

where $\sum_{l=0}^{-1}(\cdot) \equiv 0$. From (3.1e),

$$(3.3) \qquad x_t = \sum_{s=0}^{t-1} (-\beta)^s \left(\prod_{i=1}^{s} \theta_{t-i}\right) y_{t-s}$$

where $\prod_{i=1}^{0}(\cdot) \equiv 1$. Next, let $z_1 = e_1$ and, for $t \geq 2$,

$$(3.4) \qquad z_t = e_t + \theta_{t-1} \phi_{t-1}.$$

The value of the parameterization with $\beta$ is that it enables us to simultaneously obtain results for two important algorithms. When $\beta = 0$, then $x_t = y_t$ from which it follows that $\phi_t = e_t$ and $z_t = y_t$ and therefore (3.1a)–(3.1e) is PLR (AML, RML$_1$)[11, 20, 21, 24]. When $\beta = 1$, then $x_t = e_t$ and $\phi_t = e_t - \theta_{t-1} \phi_{t-1}$ and thus (3.1a)–(3.1e) is RML$_2$ [12, 21] without monitoring to ensure that each estimate $\theta_t$ is in $\Theta = (-1, 1)$.

For any $\beta$, these $\theta_t$ can be expressed in the form of a regression estimate:

$$(3.5) \qquad \theta_t = \left\{\sum_{s=2}^{t} \phi_{s-1}^2\right\}^{-1} \sum_{s=2}^{t} z_s \phi_{s-1}, \quad t \geq 2.$$

An induction argument for (3.5) goes as follows. Set $P_t = t\bar{P}_t = \sum_{s=2}^{t} \phi_{s-1}^2$. Note that from (D1), $P_t > 0$ for all $t > 1$ and therefore $P_t^{-1}$ exists a.s. From (3.1a)–(3.1e) and (3.4), $\theta_2 = \left(1/2\phi_1^2\right)^{-1} 1/2(z_2 \phi_1)$, which is (3.5) for $t = 2$. Suppose then it is true for some $t \geq 2$; i.e.,

$$(3.6) \qquad P_t \theta_t = \sum_{s=2}^{t} z_s \phi_{s-1}.$$

Then

$$P_{t+1}\theta_{t+1} = P_{t+1}(\theta_t + P_{t+1}^{-1}\phi_t e_{t+1}) = (P_t + \phi_t^2)\theta_t + \phi_t e_{t+1}$$

$$= \sum_{s=2}^{t} z_s \phi_{s-1} + \phi_t(\phi_t \theta_t + e_{t+1}) \text{ (from the induction hypothesis (3.6))}$$

$$= \sum_{s=2}^{t} z_s \phi_{s-1} + \phi_t z_{t+1} = \sum_{s=2}^{t+1} z_s \phi_{s-1}.$$

Hence, (3.5) is true for $t+1$ and by induction therefore for all $t$.

For use below, we define the stationary analogues $e_t(\theta), x_t(\theta), \phi_t(\theta)$ and $z_t(\theta)$ of $e_t, x_t, \phi_t$ and $z_t$:

(3.7) $$e_t(\theta) = (1+\theta B)^{-1} y_t,$$

(3.8) $$x_t(\theta) = (1+\theta\beta B)^{-1} y_t = \sum_{j=0}^{\infty} (-\beta\theta)^j y_{t-j},$$

$$\phi_t(\theta) = (1+\theta B)^{-1} x_t(\theta) = \sum_{j=0}^{\infty} (-\theta)^j x_{t-j}(\theta)$$

(3.9) $$= (1+\theta B)^{-1}(1+\theta\beta B)^{-1} y_t,$$

so $\phi_t(\theta) = e_t(\theta)$ when $\beta = 0$. From (3.7)–(3.9),

(3.10) $$z_t(\theta) = e_t(\theta) + \theta\phi_{t-1}(\theta) = [(1+\theta B)^{-1} + \theta B(1+\theta B)^{-1}(1+\theta\beta B)^{-1}] y_t.$$

From (3.7)–(3.10),

(3.11) $$E[\phi_t^2(\theta)] = \int_{-\pi}^{\pi} \frac{1}{|(1+\theta e^{i\omega})(1+\beta\theta e^{i\omega})|^2} g_y(\omega) d\omega,$$

(3.12) $$E[\phi_{t-1}(\theta)e_t(\theta)] = \int_{-\pi}^{\pi} \frac{e^{i\omega}}{(1+\theta e^{i\omega})(1+\beta\theta e^{i\omega})} \frac{1}{(1+\theta e^{-i\omega})} g_y(\omega) d\omega$$
$$= \int_{-\pi}^{\pi} \frac{e^{i\omega} + \beta\theta}{|(1+\theta e^{i\omega})(1+\beta\theta e^{i\omega})|^2} g_y(\omega) d\omega,$$

and

(3.13) $$E[z_t(\theta)\phi_{t-1}(\theta)] = \int_{-\pi}^{\pi} \frac{e^{i\omega} + \theta(1+\beta)}{|(1+\theta e^{i\omega})(1+\beta\theta e^{i\omega})|^2} g_y(\omega) d\omega.$$

From (1.4) and (3.12), $E[\phi_{t-1}(\theta)e_t(\theta)] = -f(\theta,\beta)$. Let $e'_t(\theta) = de_t(\theta)/d\theta$. Then, from (3.7),

(3.14) $$-e'_t(\theta) = \frac{B}{1+\theta B} e_t(\theta) = \frac{B}{(1+\theta B)^2} y_t.$$

Since

$$\frac{1}{2}\frac{d}{d\theta} E[e_t^2(\theta)] = E[e'_t(\theta) e_t(\theta)],$$

from (2.9) and (3.14), the derivative of $\bar{L}(\theta)$, $\bar{L}'(\theta)$, is obtained from

(3.15) $$-\frac{1}{2}\bar{L}'(\theta) = E[-e'_t(\theta) e_t(\theta)] = \int_{-\pi}^{\pi} \frac{e^{i\omega}}{(1+\theta e^{i\omega})^2} \frac{1}{(1+\theta e^{-i\omega})} g_y(\omega) d\omega$$
$$= \int_{-\pi}^{\pi} \frac{e^{i\omega} + \theta}{|(1+\theta e^{i\omega})^2|^2} g_y(\omega) d\omega,$$

which is (3.12) with $\beta = 1$, verifying (1.6).

As a consequence of (2.4), we note that since $|z| \leq K^* < 1$ implies $0 < 1 - K^* \leq |1 - z| \leq 1 + K^*$, for (3.11) with $|\theta| \leq K^* < 1$ we have

$$(3.16) \qquad \frac{m}{(1+K^*)^4} \leq \int_{-\pi}^{\pi} \frac{1}{|(1+\theta e^{i\omega})(1+\beta\theta e^{i\omega})|^2} g_y(\omega) d\omega \leq \frac{M}{(1-K^*)^4}.$$

## 4. The convergence theorem

The following result is a generalization of the PLR and $\text{RML}_2$ results proved in [4] for MA(1) models.

**Theorem 4.1** (Convergence theorem). *Consider a series $y_t$ for which (D1)–(D4) hold. For each $\beta$ such that $0 \leq \beta \leq 1$, assume that the recursive sequence defined by (3.1a)–(3.1e) is such that, for some random $k^* = k^*(\xi)$ and $K^* = K^*(\xi)(\xi \in \Xi)$ satisfying $0 \leq k^* < \infty$ and $0 < K^* < 1$, it holds almost surely that $|\theta_{t+k^*}| \leq K^*$ for all $t$. Then for $f(\theta, \beta)$ as in (1.4):*

(a) *The sequence $\hat{\theta}_t$ defined for $t \geq 1$ by*

$$(4.1) \qquad \hat{\theta}_t = \left[ \frac{1}{t} \sum_{s=1}^{t} \int_{-\pi}^{\pi} \frac{1}{|(1+\theta_{s+k^*}e^{i\omega})(1+\beta\theta_{s+k^*}e^{i\omega})|^2} g_y(\omega) d\omega \right]^{-1}$$
$$\times \frac{1}{t} \sum_{s=1}^{t} \int_{-\pi}^{\pi} \frac{\cos\omega + (1+\beta)\theta_{s+k^*}}{|(1+\theta_{s+k^*}e^{i\omega})(1+\beta\theta_{s+k^*}e^{i\omega})|^2} g_y(\omega) d\omega$$

*has the property that $\theta_t - \hat{\theta}_t \xrightarrow{a.s.} 0$. Hence, with probability one, there is a $t_0(\xi) \geq 1$ such that $|\hat{\theta}_t| \leq (1+K^*)/2 < 1$ holds for all $t \geq t_0(\xi)$.*

(b) *For all $t > t_0(\xi)$, $\hat{\theta}_t$ satisfies a Robbins-Monro recursion,*

$$(4.2) \qquad \hat{\theta}_t = \hat{\theta}_{t-1} - \delta_t f(\hat{\theta}_{t-1}, \beta) + \delta_t \gamma_t,$$

*with $\gamma_t \xrightarrow{a.s.} 0$, $\delta_t > 0$ a.s., $\delta_t \xrightarrow{a.s.} 0$, and $\sum_{s=t_0+1}^{\infty} \delta_s = \infty$ a.s. where $f(\theta, \beta)$ has the formula (1.4).*

(c) *From (a) and (b), it follows that, with $\Theta = (-1, 1)$, the sequence $\theta_t$ converges a.s. to the compact set*

$$(4.3) \qquad \Theta_0^\beta = \{\theta \in \Theta : f(\theta, \beta) = 0\}$$

*in the sense that, on a probability one event $\Xi_0$ that does not depend on $\beta$, for each $\xi \in \Xi_0$, the cluster points of $\theta_t(\xi)$ are contained in $\Theta_0^\beta$. Further, when $y_t$ is an invertible ARMA process, then $\Theta_0^\beta$ is finite, and $\theta(\xi) = \lim_{t \to \infty} \theta_t(\xi)$ exists for every $\xi \in \Xi_0$.*

Note from (3.5), (3.11) and (3.13) that the assertion $\theta_t - \hat{\theta}_t \xrightarrow{a.s.} 0$ in part (a) of Theorem 4.1 can be formulated as the assertion that

$$\left\{ \frac{1}{t} \sum_{s=1}^{t} \phi_{s-1}^2 \right\}^{-1} \frac{1}{t} \sum_{s=1}^{t} z_s \phi_{s-1}$$
$$- \left[ \frac{1}{t} \sum_{s=1}^{t} E[\phi_t^2(\theta_{s+k^*})] \right]^{-1} \frac{1}{t} \sum_{s=1}^{t} E[z_t(\theta_{s+k^*})\phi_{t-1}(\theta_{s+k^*})]$$

tends to zero a.s. In the expression above, $\phi_0 = 0$ and expectation is taken before evaluation at $\theta_{s+k^*}$.

The proof of Theorem 4.1, given in Section 4.2. In [5], we provide complete results concerning the existence of $k^*$ and $K^*$ with the required properties for several incorrect model examples as well as for the correct model situation for $\beta = 0$ (PLR) and provide more limited results for the case $\beta = 1$ (RML$_2$) with a particular monitoring scheme. For the latter case, we also report on simulation results which demonstrate the existence of the variates $k^*, K^*$ as in Theorem 4.1 with the consequence that monitoring becomes unnecessary for sufficiently large $t$. In the correct model case $y_t = \theta \epsilon_{t-1} + \epsilon_t$ with i.i.d. $\epsilon_t$, Lai and Ying [19] show for their monitored RML$_2$ that this happens a.s. and the conclusions of Theorem 4.1 concerning our approximating sequence (4.1) apply.

### 4.1. Preliminary results

Here we present some needed technical results. We first quote, without proof, a powerful result from martingale theory [17, Lemma 1, part (i)]. Unless specified otherwise, all limits (liminfs, limsups, etc.) are with respect to $t$ and for simplicity the $t \to \infty$ will be usually suppressed.

**Proposition 4.1.** *Let $\{\tilde{\epsilon}_t\}$ be a martingale difference sequence with respect to an increasing sequence of $\sigma$-fields $\{\mathcal{F}_t\}$ such that $\sup_t E[|\tilde{\epsilon}_t|^{2p}|\mathcal{F}_{t-1}] < \infty$ holds a.s. for some $p > 1$. Let $\tilde{z}_t$ be an $\mathcal{F}_{t-1}$-measurable random variable for every $t$. Then $\sum_{s=1}^t \tilde{z}_s \tilde{\epsilon}_s$ converges almost surely on $\{\sum_{s=1}^\infty \tilde{z}_s^2 < \infty\}$, and for every $\eta > 1/2$,*

$$\frac{\left(\sum_{s=1}^t \tilde{z}_s \tilde{\epsilon}_s\right)}{\left(\sum_{s=1}^t \tilde{z}_s^2\right)^\eta} \xrightarrow{a.s.} 0 \ \ on \ \ \left\{\sum_{s=1}^\infty \tilde{z}_s^2 = \infty\right\}.$$

Since

$$\frac{1}{t}\sum_{s=1}^t \tilde{z}_s \tilde{\epsilon}_s = \left\{\frac{\sum_{s=1}^t \tilde{z}_s \tilde{\epsilon}_s}{\sum_{s=1}^t \tilde{z}_s^2}\right\} \frac{1}{t}\sum_{s=1}^t \tilde{z}_s^2,$$

it is clear that a corollary of this Proposition is

**Proposition 4.2.** *Under the assumptions of Proposition 4.1, if $\limsup t^{-1} \times \sum_{s=1}^t \tilde{z}_s^2 < \infty$ a.s., then $t^{-1} \sum_{s=1}^t \tilde{z}_s \tilde{\epsilon}_s \xrightarrow{a.s.} 0$.*

Recall from (2.1) that $y_t = \epsilon_t + \sum_{s=1}^\infty \kappa_s \epsilon_{t-s}$ since $\kappa_0 = 1$. A second consequence of Proposition 4.1 is

**Proposition 4.3.** *Suppose that the m.d.s. $\epsilon_t$ in (D2) is such that $\sup_t E[|\epsilon_t|^{2p}|\mathcal{F}_{t-1}] < \infty$ holds a.s. for some $p > 1$. Then for any sequence $\hat{y}_t = y_t - \tilde{y}_{t-1}$ in which $\tilde{y}_{t-1}$ is $\mathcal{F}_{t-1}$-measurable, it holds that $\liminf t^{-1} \sum_{s=1}^t \hat{y}_s^2 \geq \sigma_\epsilon^2$ a.s., where $\sigma_\epsilon^2 = E[\epsilon_t^2]$.*

*Proof.* From (2.1), $\hat{y}_t = y_t - \tilde{y}_{t-1} = \epsilon_t + \tilde{z}_t$ where $\tilde{z}_t = -\tilde{y}_{t-1} + \sum_{s=1}^\infty \kappa_s \epsilon_{t-s}$ is $\mathcal{F}_{t-1}$-measurable since $\sum_{s=1}^\infty \kappa_s \epsilon_{t-s}$ is $\mathcal{F}_{t-1}$-measurable by (2.2) and $\tilde{y}_{t-1}$ is $\mathcal{F}_{t-1}$-measurable by assumption. Then

(4.4)
$$\begin{aligned}\frac{1}{t}\sum_{s=1}^t \hat{y}_s^2 &= \frac{1}{t}\sum_{s=1}^t \epsilon_s^2 + \frac{2}{t}\sum_{s=1}^t \epsilon_s \tilde{z}_s + \frac{1}{t}\sum_{s=1}^t \tilde{z}_s^2 \\ &= \frac{1}{t}\sum_{s=1}^t \epsilon_s^2 + \left\{2\frac{\sum_{s=1}^t \epsilon_s \tilde{z}_s}{\sum_{s=1}^t \tilde{z}_s^2} + 1\right\}\frac{1}{t}\sum_{s=1}^t \tilde{z}_s^2.\end{aligned}$$

Consider first the event that $\sum_{s=1}^{t} \tilde{z}_s^2 \xrightarrow{a.s.} l < \infty$. Then $t^{-1} \sum_{s=1}^{t} \tilde{z}_s^2 \xrightarrow{a.s.} 0$ and, by the preceding Proposition, $t^{-1} \sum_{s=1}^{t} \epsilon_s \tilde{z}_s \xrightarrow{a.s.} 0$. Hence, from (2.5) and the first equation in (4.4), $\lim t^{-1} \sum_{s=1}^{t} \hat{y}_s^2 = t^{-1} \sum_{s=1}^{t} \epsilon_s^2 = \sigma_\epsilon^2$ so the assertion holds in this event. In the complementary event, $\sum_{s=1}^{t} \tilde{z}_s^2 \xrightarrow{a.s.} \infty$, from (4.4), it follows that

$$
\begin{aligned}
(4.5) \quad \liminf \frac{1}{t} \sum_{s=1}^{t} \hat{y}_s^2 &= \liminf \left( \frac{1}{t} \sum_{s=1}^{t} \epsilon_s^2 + \left\{ 2 \frac{\sum_{s=1}^{t} \epsilon_s \tilde{z}_s}{\sum_{s=1}^{t} \tilde{z}_s^2} + 1 \right\} \frac{1}{t} \sum_{s=1}^{t} \tilde{z}_s^2 \right) \\
&= \sigma_\epsilon^2 + \liminf \left( \left\{ 2 \frac{\sum_{s=1}^{t} \epsilon_s \tilde{z}_s}{\sum_{s=1}^{t} \tilde{z}_s^2} + 1 \right\} \frac{1}{t} \sum_{s=1}^{t} \tilde{z}_s^2 \right) \quad a.s.
\end{aligned}
$$

By Proposition 4.1, $\sum_{s=1}^{t} \epsilon_s \tilde{z}_s / \sum_{s=1}^{t} \tilde{z}_s^2 \xrightarrow{a.s.} 0$. Hence, the second expression in (4.5) is nonnegative, and the proof is complete. □

**Proposition 4.4.** *Under (2.4), for each $\beta \in [0,1]$, the function $f(\theta, \beta)$ defined by (1.4) is infinitely differentiable on $\Theta = (-1, 1)$, and $\Theta_0^\beta$ defined by (4.3) is a nonempty compact subset of $\Theta$. In the case $\beta = 1$, $\Theta_0^1$ contains the (nonempty) set of minimizers over $\Theta$ of $\bar{L}(\theta)$ defined by (2.9).*

*Proof.* The differentiability assertion follows from (2.4) via the dominated convergence theorem. Except for compactness of $\Theta_0^1$, which will be discussed below, the assertions concerning $\bar{L}(\theta)$ and $f(\theta, 1)$ were obtained subsequent to (2.10). The remaining assertions follow from the continuity of $f(\theta, \beta)$ and the limit properties

$$
(4.6) \quad \lim_{\theta \to -1} f(\theta, \beta) = -\infty
$$

and

$$
(4.7) \quad \lim_{\theta \to 1} f(\theta, \beta) = \infty.
$$

Indeed, from (4.6)–(4.7), for any $K > 0$ there exists an $0 < \epsilon(K, \beta) < 1$ such that $f(\theta, \beta) \leq -K$ for all $\theta \in (-1, -1+\epsilon)$ and $f(\theta, \beta) \geq K$ for all $\theta \in (1-\epsilon, 1)$. Therefore $f(\theta, \beta)$ must change sign over $[-1+\varepsilon, 1-\varepsilon]$. Hence $f(\theta, \beta)$ is non-constant and has a zero in this interval and, moreover, $\Theta_0^\beta \subseteq [-1+\varepsilon, 1-\varepsilon]$. Finally, since $f(\theta, \beta)$ is continuous on this interval, $\Theta_0^\beta$ is compact. An analogous argument applies to $\Theta_0^1$.

To verify (4.6), we note that $g_y(\omega) = g_y(-\omega), -\pi \leq \omega \leq \pi$ yields

$$
f(\theta, \beta) = -\int_{-\pi}^{\pi} \frac{\cos \omega + \beta \theta}{|(1 + \theta e^{i\omega})(1 + \beta \theta e^{i\omega})|^2} g_y(\omega) d\omega.
$$

Because $0 \leq \beta < 1$, for $0 < \varepsilon < 1 - \beta$ there is a $\delta = \delta(\varepsilon) \in (0, \pi)$ such that $\cos \omega + \beta \theta \geq \varepsilon$ whenever $|\omega| \leq \delta$ and $-1 \leq \theta \leq 0$. For such $\varepsilon, \delta$, we obtain

$$
\begin{aligned}
&\lim_{\theta \to -1} \int_{-\pi}^{\pi} \frac{\cos \omega + \beta \theta}{|(1 + \theta e^{i\omega})(1 + \beta \theta e^{i\omega})|^2} g_y(\omega) d\omega \\
(4.8) \quad &= \left\{ \int_{-\pi}^{-\delta} + \int_{\delta}^{\pi} \right\} \frac{\cos \omega + \beta \theta}{|(1 - e^{i\omega})(1 - \beta e^{i\omega})|^2} g_y(\omega) d\omega \\
(4.9) \quad &\quad + \lim_{\theta \to -1} \int_{-\delta}^{\delta} \frac{\cos \omega + \beta \theta}{|(1 + \theta e^{i\omega})(1 + \beta \theta e^{i\omega})|^2} g_y(\omega) d\omega \\
&= \infty,
\end{aligned}
$$

because (4.8) is finite, whereas for (4.9) we have

$$\lim_{\theta \to -1} \int_{-\delta}^{\delta} \frac{\cos \omega + \beta \theta}{\left|(1 + \theta e^{i\omega})(1 + \beta \theta e^{i\omega})\right|^2} g_y(\omega) \, d\omega$$
$$\geq \varepsilon \, m \lim_{\theta \to -1} \int_{-\delta}^{\delta} \left|(1 + \theta e^{i\omega})(1 + \beta \theta e^{i\omega})\right|^{-2} d\omega = \infty.$$

This yields (4.6), and (4.7) follows by an analogous argument. □

**Proposition 4.5.** *Let $y_t$ be an invertible ARMA process, then for each $\beta \in [0, 1]$, the set $\Theta_0^\beta = \{\theta \in (-1, 1) : f(\theta, \beta) = 0\}$ is finite.*

*Proof.* $\kappa(z)$ in (D2) has the form $\kappa(z) = \eta(z)/\phi(z)$ where $\eta(z)$ and $\phi(z)$ are polynomials, of degrees $d_\eta$ and $d_\phi$, respectively, having no common zeroes and having all zeros in $\{|z| > 1\}$. Setting $z = e^{i\omega}$ and $h(z) = (1 + \theta z)(1 + \beta \theta z)$, we obtain from $dz = iz\,d\omega$ that

$$-f(\theta, \beta) = \int_{-\pi}^{\pi} \frac{e^{i\omega} + \beta \theta}{\left|(1 + \theta e^{i\omega})(1 + \beta \theta e^{i\omega})\right|^2} g_y(\omega) \, d\omega$$
$$= \frac{\sigma_\varepsilon^2}{2\pi i} \int_{|z|=1} \frac{(z + \beta\theta) \eta(z) \eta(z^{-1})}{zh(z) h(z^{-1}) \phi(z) \phi(z^{-1})} dz$$
$$= \frac{\sigma_\varepsilon^2}{2\pi i} \int_{|z|=1} z^{1+d_\phi - d_\eta} \frac{(z + \beta\theta) \eta(z) \{z^{d_\eta} \eta(z^{-1})\}}{h(z) \{z^2 h(z^{-1})\} \phi(z) \{z^{d_\phi} \phi(z^{-1})\}} dz.$$

The function

$$w(z) = \sigma_\varepsilon^2 z^{1+d_\phi - d_\eta} \frac{(z + \beta\theta) \eta(z) \{z^{d_\eta} \eta(z^{-1})\}}{h(z) \{z^2 h(z^{-1})\} \phi(z) \{z^{d_\phi} \phi(z^{-1})\}}$$

is nonzero on $\{|z| = 1\}$ and has poles interior to the unit circle at $-\theta$, $-\beta\theta$, at the zeroes of $z^{d_\phi} \phi(z^{-1})$, and, if $1 + d_\phi - d_\eta < 0$, also at 0. If $z_j, j = 1, \ldots, n$ are the distinct poles in $\{z : |z| < 1\}$, then, by the Residue Theorem of complex analysis, e.g., (4.7-10) of Henrici [13], it follows that

$$f(\theta, \beta) = -\sum_{j=1}^{n} \operatorname{Res}_{z=z_j} w(z),$$

where, if $z_j$ is a pole of order $J \geq 1$,

$$\operatorname{Res}_{z=z_j} w(z) = \frac{1}{(J-1)!} \lim_{z \to z_j} \frac{d^{J-1}}{dz^{J-1}} \left\{(z - z_j)^J w(z)\right\}.$$

Thus each $\operatorname{Res}_{z=z_j} w(z)$ is a rational function of $\theta$, and therefore the same is true of $f(\theta, \beta)$. Consequently, $f(\theta, \beta) = 0$ holds for only finitely many $\theta$ in $(-1, 1)$. □

The final preliminary result addresses convergence of a Robbins-Monro type recursion that will be applied to demonstrate convergence of the general recursive algorithm. It is a special case of a correction and extension by Findley [9] of a result that is implicit in the proof of a theorem of Fradkov presented in Derevitzkiĭ and Fradkov [8] for the case of monotonically decreasing $\delta_t$. The result below is also implicit in the proofs of Theorem 2.2.2 and Corollary 2.2.1 of Chen [7] which cover the case of vector $\theta$ more completely than Findley [9].

**Proposition 4.6.** *Let $\hat{\theta}_t$, $t \geq t_0$ be a non-stochastic, real-valued sequence satisfying*

$$\hat{\theta}_t = \hat{\theta}_{t-1} - \delta_t f(\hat{\theta}_{t-1}) + \delta_t \gamma_t, \ t > t_0$$

*for some real-valued function $f(\theta)$, with $\gamma_t$, $t > t_0$ satisfying $\gamma_t \to 0$ and with $\delta_t$, $t \geq t_0$ satisfying $\delta_t \geq 0$, $\delta_t \to 0$, and $\sum_{t=t_0+1}^{\infty} \delta_t = \infty$. Suppose there is a bounded open set $\tilde{\Theta}$ on which $f(\theta)$ is continuously differentiable and which is such that the sequence $\hat{\theta}_t$ enters $\tilde{\Theta}$ infinitely often and has no cluster point on the boundary of $\tilde{\Theta}$. Then $\hat{\theta}_t$ is bounded, and its cluster points belong to $\tilde{\Theta}_0 = \{\theta \in \tilde{\Theta}: f(\theta) = 0\}$, i.e., $\hat{\theta}_t \to \tilde{\Theta}_0$. The set of cluster points is compact. If $\tilde{\Theta}_0$ is finite, then $\hat{\theta}_t$ converges to some $\theta \in \tilde{\Theta}_0$.*

### 4.2. Proof of the convergence theorem

The proof of Theorem 4.1 follows from a set of technical lemmas and propositions given below. Proposition 4.7 provides a set of technical results needed to prove the Theorem's two main assertions: (i) the asymptotic equivalence of $\theta_t$ and the sequence $\hat{\theta}_t$ (Proposition 4.8) and (ii) (Proposition 4.9) the fact that $\hat{\theta}_t$ satisfies a.s. a Robbins-Monro recursion of the form considered in Proposition 4.6.

Hereafter, $K$ or sometimes $k$ (or these letters with decorations) will denote a generic upper bound (not always the same one) that is finite, or when it is random, finite a.s. A random $K$ will be shown as $K(\xi)$ with $\xi \in \Xi$ on first appearance whenever the randomness is not immediately clear from context. Again, unless specified otherwise, all limits (liminfs, limsups, etc.) are with respect to $t$ and usually the $t \to \infty$ will be omitted. The notation $o_{a.s.}(1)$ denotes convergence to zero with probability one.

**Proposition 4.7.** *Under the assumptions of Theorem 4.1, for the general recursive algorithm, the assertions* (a)–(c) *below follow:*

(a) $\liminf t^{-1} \sum_{s=1}^{t} \phi_s^2 \geq \sigma_\epsilon^2$ *a.s. and* $(t^{-1} \sum_{s=1}^{t} \phi_s^2)^{-1} \leq K(\xi) < \infty$, *and thus, from (3.1b), $\bar{P}_t^{-1}$ is bounded a.s.*

(b) *For $t \geq 1$, $e_t = \sum_{j=0}^{\infty} \kappa_j^e(t) \epsilon_{t-j}$; $\phi_t = \sum_{j=0}^{\infty} \kappa_j^\phi(t) \epsilon_{t-j}$; $x_t = \sum_{j=0}^{\infty} \kappa_j^x(t) \epsilon_{t-j}$; and $z_t = \sum_{j=0}^{\infty} \kappa_j^z(t) \epsilon_{t-j}$ where for every $j$, $\kappa_j^e(t), \kappa_j^\phi(t), \kappa_j^x(t)$ and $\kappa_j^z(t)$ are $\mathcal{F}_{t-1}$-measurable. Moreover, there exist $\tilde{\kappa}_j$ such that*

$$\max_j \{|\kappa_j^e(t)|, |\kappa_j^\phi(t)|, |\kappa_j^x(t)|, |\kappa_j^z(t)|\} \leq \tilde{\kappa}_j$$

*and $\sum_j^{\infty} \tilde{\kappa}_j < \infty$ a.s. Hence, the sequences $e_t, \phi_t, x_t$ and $z_t$ are uniformly bounded a.s.*

(c) $\theta_t - \theta_{t-1} = o_{a.s.}(1)$.

*Proof of* (a). From (3.1d), $\phi_t = x_t - \theta_{t-1} e_{t-1} = y_t - \theta_{t-1}(\beta x_{t-1} + e_{t-1})$. Since $\theta_{t-1}(\beta x_{t-1} + e_{t-1})$ is $\mathcal{F}_{t-1}$-measurable, by Proposition 4.3,

$$(4.10) \qquad \liminf t^{-1} \sum_{s=1}^{t} \phi_s^2 \geq \sigma_\epsilon^2 \ a.s.$$

Continuing, from (4.10), for any $0 < L_1 < \sigma_\epsilon^2$, there exists $t_0 = t_0(L_1, \xi)$ such that $t^{-1} \sum_{s=1}^{t} \phi_s^2 > L_1$ a.s. for all $t \geq t_0$. Let $L_2(\xi) \equiv \min_{1 \leq t < t_0} t^{-1} \sum_{s=1}^{t} \phi_s^2$. Then $0 < L_2 < \infty$ a.s. This follows since $t_0$ is finite and $\phi_t$ is a finite valued sequence with probability one, hence $L_2 < \infty$. Moreover, since $\phi_1 = y_1$, under (D1) it follows that $L_2 > 0$ a.s. Hence, $(t^{-1} \sum_{s=1}^{t} \phi_s^2)^{-1} \leq \max\{L_1^{-1}, L_2^{-1}\} < \infty$ a.s. and the proof of part (a) is complete. □

*Proof of* (b). Set $\theta_0 = 0$. From $e_1 = y_1$ and $e_t = y_t - \theta_{t-1}e_{t-1}, t \geq 2$, it follows that $\kappa_j^e(1) = \kappa_j$ for all $j$, that $\kappa_0^e(t) = \kappa_0$ for all $t \geq 1$, and that $\kappa_j^e(t) = \kappa_j(t) - \theta_{t-1}\kappa_j^e(t-1)$ for all $t \geq 2, j \geq 1$. It follows by induction that

$$(4.11) \qquad \kappa_j^e(t) = \sum_{l=0}^{\min(j,t-1)} (-1)^l \kappa_{j-l} \prod_{i=1}^{l} \theta_{t-i} \quad \text{where } \prod_{i=1}^{0}(\cdot) \equiv 1.$$

Since for some $k^*$ finite, $|\theta_{t+k^*}| < 1$ for all $t \geq 1$, we have that $|\theta_t| \leq K(\xi) < \infty$. First suppose that $K < 1$. Then from (4.11),

$$|\kappa_j^e(t)| \leq \sum_{l=0}^{\min(j,t-1)} |\kappa_{j-l}| \prod_{i=1}^{l} |\theta_{t-i}| \leq \sum_{l=0}^{j} K^l |\kappa_{j-l}|$$

and since $K < 1$, $\sum_{j=0}^{\infty} |\kappa_j^e(t)| \leq \sum_{j=0}^{\infty} \sum_{l=0}^{j} K^l |\kappa_{j-l}| = \sum_{l=0}^{\infty} K^l \sum_{p=0}^{\infty} |\kappa_p| < \infty$ where $p = j - l$. So the result holds for the case of $0 < K < 1$.

Otherwise, suppose $1 \leq K < \infty$. For all $t \geq k^*$, we have that $|\theta_t| \leq K^*(\xi) < 1$, so $K(\xi) = \lambda(\xi)K^*(\xi)$ for $\lambda > 1$. For simplicity of notation, replace $K^*$ by $\rho$. We next show that $\prod_{i=1}^{l} |\theta_{t-i}| \leq \lambda^{k^*}\rho^l$ for $l \leq t$. First suppose $t \leq k^*$. Then $\prod_{i=1}^{l} |\theta_{t-i}| \leq \lambda^l \rho^l \leq \lambda^{k^*}\rho^l$. Next suppose $t > k^*$ and $l \leq t - k^*$. Then, $\prod_{i=1}^{l} |\theta_{t-i}| \leq \rho^l < \rho^l \lambda^{k^*}$ since $|\theta_{t-i}| \leq \rho$ for $1 \leq i \leq t - s^*$. Finally, suppose $t > k^*$ and $l > t - s^*$. Then since $l \leq t$,

$$\prod_{i=1}^{l} |\theta_{t-i}| = \prod_{i=1}^{t-s^*} |\theta_{t-i}| \prod_{i=t-s^*+1}^{l} |\theta_{t-i}| \leq \rho^{t-s^*}\lambda^{l-(t-s^*)}\rho^{l-(t-s^*)}$$
$$= \rho^l \lambda^{l-(t-s^*)} = \lambda^{k^*}\lambda^{l-t}\rho^l \leq \lambda^{k^*}\rho^l.$$

Hence, generally $\prod_{i=1}^{l} |\theta_{t-i}| \leq \lambda^{k^*}\rho^l$. Setting $\kappa_j^e(\xi) = \lambda^{k^*} \sum_{l=0}^{j} \rho^l |\kappa_{j-l}|$, we have

$$|\kappa_j^e(t)| \leq \sum_{l=0}^{j} |\kappa_{j-l}| \prod_{i=1}^{l} |\theta_{t-i}| \leq \lambda^{k^*} \sum_{l=0}^{j} \rho^l |\kappa_{j-l}| = \kappa_j^e,$$

and since $|\rho| < 1$, $\sum_{j=0}^{\infty} \kappa_j^e < \infty$ a.s.

Next, from (3.3)

$$(4.12) \qquad \kappa_j^x(t) = \sum_{l=0}^{\min(j,t-1)} (-\beta)^l \kappa_{j-l} \prod_{i=1}^{l} \theta_{t-i},$$

and since $0 \leq \beta \leq 1$, an argument like that for $e_t$ can be applied and to obtain the existence of a $\kappa_j^x$ such that

$$(4.13) \qquad |\kappa_j^x(t)| \leq \kappa_j^x \text{ and } \sum_{j=0}^{\infty} \kappa_j^x < \infty \, a.s.$$

Continuing, since $\phi_1 = x_1$ and $\phi_t = x_t - \theta_{t-1}\phi_{t-1}$ for $t \geq 2$, it follows similarly that

$$(4.14) \qquad \kappa_j^\phi(t) = \sum_{l=0}^{\min(j,t-1)} (-1)^l \kappa_{j-l}^x(t) \prod_{i=1}^l \theta_{t-i}.$$

From (4.12) and (4.13), substituting $\kappa_j^x(t)$ for $\kappa_j$, the same kind of argument can be applied to (4.14) to yield

$$(4.15) \qquad |\kappa_j^\phi(t)| \leq \kappa_j^\phi \text{ with } \sum_{j=0}^\infty \kappa_j^\phi < \infty \text{ a.s.}$$

Finally, for $t \geq 2$, we have, from $z_t = e_t + \theta_{t-1}\phi_{t-1}$,

$$\sum_{j=0}^\infty \kappa_j^z(t)\epsilon_{t-j} = \sum_{j=0}^\infty \kappa_j^e(t)\epsilon_{t-j} + \theta_{t-1}\sum_{j=0}^\infty \kappa_j^\phi(t-1)\epsilon_{t-1-j},$$

for $t \geq 2$ from which it follows that

$$(4.16) \qquad \kappa_j^z(t) = \kappa_j^e(t) + \theta_{t-1}\kappa_{j-1}^\phi(t-1),$$

where $\kappa_{-1}^\phi(t) \equiv 0$. Since $\sup_t |\theta_t| < \infty$ a.s.,

$$|\kappa_j^z(t)| \leq \kappa_j^e + \sup_t |\theta_t| \kappa_{j-1}^\phi \quad a.s.,$$

where $\kappa_{-1}^\phi \equiv 0$, so there is a $\kappa_j^z$ such that $|\kappa_j^z(t)| \leq \kappa_j^z$ and $\sum_{j=0}^\infty |\kappa_j^z| < \infty$ a.s. for $t \geq 2$. Since $z_1 = e_1$, it thus follows that $\tilde{\kappa}_j = \max_j\{|\kappa_j^e|, |\kappa_j^\phi|, |\kappa_j^x|, |\kappa_j^z|\}$ satisfies $\sum_j^\infty \tilde{\kappa}_j < \infty$ a.s.

From this, we see that $e_t, \phi_t, x_t$ and $z_t$ are bounded a.s. For example,

$$|\phi_t| = \left|\sum_{j=0}^\infty \kappa_j^\phi(t)\epsilon_{t-j}\right| \leq \sup_{-\infty < t < \infty} |\epsilon_t| \sum_{j=0}^\infty \tilde{\kappa}_j < \infty \quad a.s.$$

From (4.11)–(4.12), (4.14) and (4.16), $\kappa_j^e(t), \kappa_j^\phi(t), \kappa_j^x(t)$ and $\kappa_j^z(t)$ are each $\mathcal{F}_{t-1}$-measurable for every $j$. Hence, part (b) of the Proposition is proved. $\square$

*Proof of* (c). By parts (a) and (b), $|\theta_t - \theta_{t-1}| \leq t^{-1}\bar{P}_t^{-1}|e_t||\phi_{t-1}| \leq t^{-1}K(\xi)$ where $K(\xi) < \infty$ and thus part (c) follows and the proof of Proposition 4.7 is complete. $\square$

**Lemma 4.1.** *Under the assumptions of Theorem 4.1, we have:*

(a) *If $\tilde{\kappa}_j(t)$ are $\mathcal{F}_{t-1}$-measurable such that $|\tilde{\kappa}_j(t)| \leq \tilde{\kappa}_j$ for $j \geq 0$, with $\sum_{j=0}^\infty \tilde{\kappa}_j < \infty$ a.s., then for all $p \geq 1$ and each $0 \leq j < \infty$,*

$$(4.17) \qquad \frac{1}{t}\sum_{s=2}^t \left(\sum_{l=1}^{\min(j,s-1)} \tilde{\kappa}_{j-l}(s) \prod_{i=1}^l (s-i)^{-1}\bar{P}_{s-i}^{-1}\phi_{s-i-1}e_{s-i}\right)^p \xrightarrow{a.s.} 0,$$

*and*

$$(4.18) \qquad \frac{1}{t}\sum_{s=2}^t \left(\sum_{l=1}^{\min(j,s-1)} \tilde{\kappa}_{j-l}(s) \sum_{i=0}^l (s-i)^{-1}\bar{P}_{s-i}^{-1}\phi_{s-i-1}e_{s-i}\right)^p \xrightarrow{a.s.} 0.$$

*In particular,*

$$(4.19) \quad \frac{1}{t}\sum_{s=2}^{t}\left(\sum_{l=1}^{\min(j,s-1)}\tilde{\kappa}_{j-l}(s)\prod_{i=1}^{l}(s-i)^{-1}\bar{P}_{s-i}^{-1}\phi_{s-i-1}e_{s-i}\right)^{p}\epsilon_{s-j}^{2} \xrightarrow{a.s.} 0.$$

(b) *For any $0 \leq j < \infty$ and $i \leq j$,*

$$(4.20) \quad \frac{1}{t}\sum_{s=1}^{t}(\kappa_{j}^{\phi}(s))^{2}\epsilon_{s-j}^{2} = \frac{1}{t}\sum_{s=i+1}^{t}(\kappa_{j}^{\phi}(s-i))^{2}\epsilon_{s-j}^{2} + o_{a.s.}(1).$$

(c) *For $0 \leq j, l < \infty$ and $j \neq l$, then*

$$(4.21) \quad \frac{1}{t}\sum_{s=\max(j+2,l+2)}^{t}\kappa_{j}^{\phi}(s)\epsilon_{s-j}\kappa_{l}^{\phi}(s)\epsilon_{s-l} \xrightarrow{a.s.} 0.$$

*Proof of* (a). By the boundedness of $\bar{P}_{t}^{-1}, \phi_{t}, e_{t}$ (Proposition 4.7) and since $|\tilde{\kappa}_{m}(t)| \leq \tilde{\kappa}_{m}$ for all $m \geq 0$ and $t \geq 1$,

$$\frac{1}{t}\sum_{s=2}^{t}\left(\sum_{l=1}^{\min(j,s-1)}\tilde{\kappa}_{j-l}(s)\prod_{i=1}^{l}(s-i)^{-1}\bar{P}_{s-i}^{-1}\phi_{s-i-1}e_{s-i}\right)^{p}$$

$$\leq \frac{1}{t}\sum_{s=2}^{t}\left(\sum_{m=0}^{j}\tilde{\kappa}_{j}\sum_{l=1}^{\min(j,s-1)}\prod_{i=1}^{l}(s-i)^{-1}|\bar{P}_{s-i}^{-1}|\|\phi_{s-i-1}\|\|e_{s-i}\|\right)^{p}$$

$$\leq K(\xi)\frac{1}{t}\sum_{s=2}^{t}\left(\sum_{l=1}^{\min(j,s-1)}\prod_{i=1}^{l}(s-i)^{-1}\right)^{p}.$$

And since for all $j \geq 0$, $p \geq 1$,

$$\frac{1}{t}\sum_{s=2}^{t}\left(\sum_{l=1}^{\min(j,s-1)}\prod_{i=1}^{l}(s-i)^{-1}\right)^{p} \leq \frac{K}{t}\sum_{s=2}^{t}(s-\min(j,s-1))^{-p} \longrightarrow 0,$$

(4.17) follows, as does (4.19), by the boundedness of $\epsilon_{t}$. Similarly,

$$\frac{1}{t}\sum_{s=2}^{t}\left(\sum_{l=1}^{\min(j,s-1)}\tilde{\kappa}_{j-l}(s)\sum_{i=0}^{l}(s-i)^{-1}\bar{P}_{s-i}^{-1}\phi_{s-i-1}e_{s-i}\right)^{p}$$

$$(4.22) \quad \leq K(\xi)\frac{1}{t}\sum_{s=2}^{t}\left(\sum_{l=1}^{\min(j,s-1)}\sum_{i=0}^{l}(s-i)^{-1}\right)^{p}$$

$$\leq K(\xi)\frac{K}{t}\sum_{s=2}^{t}\left(\sum_{i=0}^{\min(j,s-1)}(s-i)^{-1}\right)^{p} \longrightarrow 0,$$

and (4.18) follows. □

*Proof of* (b). From (4.12) and the recursion (3.1a) for $\theta_t$, we have, for $s \geq j+2$,

$$
\begin{aligned}
\kappa_j^x(s) &= \sum_{l=0}^{j}(-\beta)^l \kappa_{j-l} \prod_{i=1}^{l} \theta_{s-i} \\
&= \sum_{l=0}^{j}(-\beta)^l \kappa_{j-l} \prod_{i=1}^{l} \left(\theta_{s-i-1} + (s-i)^{-1}\bar{P}_{s-i}^{-1}\phi_{s-i-1}e_{s-i}\right) \\
&= \sum_{l=0}^{j}(-\beta)^l \kappa_{j-l} \prod_{i=1}^{l} \theta_{s-i-1} + \sum_{l=0}^{j}(-\beta)^l \kappa_{j-l} \prod_{i=1}^{l}(s-i)^{-1}\bar{P}_{s-i}^{-1}\phi_{s-i-1}e_{s-i} \\
&= \kappa_j^x(s-1) + w_j^x(s).
\end{aligned}
\tag{4.23}
$$

where

$$
w_j^x(s) = \sum_{l=0}^{j}(-\beta)^l \kappa_{j-l} \prod_{i=1}^{l}(s-i)^{-1}\bar{P}_{s-i}^{-1}\phi_{s-i-1}e_{s-i},
\tag{4.24}
$$

Continuing, from (4.14) and (4.23)–(4.24), for $s \geq j+2$,

$$
\begin{aligned}
\kappa_j^\phi(s) &= \sum_{l=0}^{j}(-1)^l \kappa_{j-l}^x(s) \prod_{i=1}^{l} \theta_{s-i} \\
&= \sum_{l=0}^{j}(-1)^l \kappa_{j-l}^x(s) \prod_{i=1}^{l} \left(\theta_{s-i-1} + (s-i)^{-1}\bar{P}_{s-i}^{-1}\phi_{s-i-1}e_{s-i}\right) \\
&= \sum_{l=0}^{j}(-1)^l \kappa_{j-l}^x(s) \prod_{i=1}^{l} \theta_{s-i-1} \\
&\quad + \sum_{l=0}^{j}(-1)^l \kappa_{j-l}^x(s) \prod_{i=1}^{l}(s-i)^{-1}\bar{P}_{s-i}^{-1}\phi_{s-i-1}e_{s-i} \\
&= \sum_{l=0}^{j}(-1)^l (\kappa_{j-l}^x(s-1) + w_{j-l}^x(s)) \prod_{i=1}^{l} \theta_{s-i-1} \\
&\quad + \sum_{l=0}^{j}(-1)^l \kappa_{j-l}^x(s) \prod_{i=1}^{l}(s-i)^{-1}\bar{P}_{s-i}^{-1}\phi_{s-i-1}e_{s-i} \\
&= \kappa_j^\phi(s-1) + w_j^\phi(s),
\end{aligned}
\tag{4.25}
$$

where from (4.24),

$$
\begin{aligned}
w_j^\phi(s) &= \sum_{l=0}^{j}(-1)^l w_{j-l}^x(s) \prod_{i=1}^{l} \theta_{s-i-1} \\
&\quad + \sum_{l=0}^{j}(-1)^l \kappa_{j-l}^x(s) \prod_{i=1}^{l}(s-i)^{-1}\bar{P}_{s-i}^{-1}\phi_{s-i-1}e_{s-i} \\
&= \sum_{l=0}^{j}(-1)^l \sum_{m=0}^{j-l}(-\beta)^m \kappa_{j-l-m} \prod_{i=1}^{m}(s-i)^{-1}\bar{P}_{s-i}^{-1}\phi_{s-i-1}e_{s-i} \prod_{n=1}^{l} \theta_{s-n-1} \\
&\quad + \sum_{l=0}^{j}(-1)^l \kappa_{j-l}^x(s) \prod_{i=1}^{l}(s-i)^{-1}\bar{P}_{s-i}^{-1}\phi_{s-i-1}e_{s-i}.
\end{aligned}
\tag{4.26}
$$

By (4.19) and (4.25)–(4.26),

$$\frac{1}{t}\sum_{s=j+2}^{t}(\kappa_j^\phi(s))^2\epsilon_{s-j}^2 = \frac{1}{t}\sum_{s=j+2}^{t}\left((\kappa_j^\phi(s-1))^2 + 2\kappa_j^\phi(s-1)w_j^\phi(s) + (w_j^\phi(s))^2\right)\epsilon_{s-j}^2.$$

Applying an argument similar to that used for part (a), it follows by the boundedness of $\beta$ and $\theta_t$ and the Cauchy-Schwarz inequality that $t^{-1}\sum_{s=j+2}^{t}(2\kappa_j^\phi(s-1)w_j^\phi(s) + (w_j^\phi(s))^2)\epsilon_{s-j}^2 = o_{a.s.}(1)$. Hence,

$$\frac{1}{t}\sum_{s=j+2}^{t}(\kappa_j^\phi(s))^2\epsilon_{s-j}^2 = \frac{1}{t}\sum_{s=j+2}^{t}(\kappa_j^\phi(s-1))^2\epsilon_{s-j}^2 + o_{a.s.}(1).$$

Finally, since $j$ is finite, then for $i \leq j$, it follows by applying the recursion (4.25) in $\kappa_j^\phi(t)$ $i-1$ additional times that (4.20) holds, because a finite sum of $o_{a.s.}(1)$ terms is $o_{a.s.}(1)$. □

*Proof of* (c). By parts (a) and (b), for $j \neq l$,

$$\frac{1}{t}\sum_{s=\max(j+2,l+2)}^{t}\kappa_j^\phi(s)\epsilon_{s-j}\kappa_l^\phi(s)\epsilon_{s-l}$$

$$= \frac{1}{t}\sum_{s=\max(j+2,l+2)}^{t}\left\{\left(\kappa_j^\phi(s-1) + \sum_{p=0}^{\min(j,s-1)}(-1)^p\kappa_{j-p}^x(s)\right.\right.$$

$$\left.\times\prod_{q=1}^{l}(s-q)^{-1}\bar{P}_{s-q}^{-1}\phi_{s-q-1}e_{s-q}\right)$$

(4.27)

$$\times\left(\kappa_l^\phi(s-1) + \sum_{r=0}^{\min(j,s-1)}(-1)^r\kappa_{j-r}^x(s)\right.$$

$$\left.\left.\times\prod_{m=1}^{r}(s-m)^{-1}\bar{P}_{s-m}^{-1}\phi_{s-m-1}e_{s-m}\right)\epsilon_{s-j}\epsilon_{s-l}\right\}$$

$$= \frac{1}{t}\sum_{s=\max(j+2,l+2)}^{t}\kappa_j^\phi(s-1)\kappa_l^\phi(s-1)\epsilon_{s-j}\epsilon_{s-l} + o_{a.s.}(1).$$

Without loss of generality, suppose $j < l < \infty$. From parts (a)–(b) and applying the argument that led to (4.27) $j-1$ additional times, we have that

$$\frac{1}{t}\sum_{s=1}^{t}\kappa_j^\phi(s)\kappa_l^\phi(s)\epsilon_{s-j}\epsilon_{s-l} = \frac{1}{t}\sum_{s=j+1}^{t}\kappa_j^\phi(s-j)\kappa_l^\phi(s-j)\epsilon_{s-j}\epsilon_{s-l} + o_{a.s.}(1)$$

$$= \frac{1}{t}\sum_{s=1}^{t-j}\kappa_j^\phi(s)\kappa_l^\phi(s)\epsilon_{s-(l-j)}\epsilon_s + o_{a.s.}(1),$$

$$= \frac{1}{t}\sum_{s=1}^{t}\kappa_j^\phi(s)\kappa_l^\phi(s)\epsilon_{s-(l-j)}\epsilon_s + o_{a.s.}(1),$$

since by (D4) and the fact that $|\kappa_m^\phi(t)| \leq K(\xi) < \infty$ for all $m \geq 0$,

$$t^{-1}\sum_{s=t-j+1}^{t}\kappa_j^\phi(s)\kappa_l^\phi(s)\epsilon_{s-(l-j)}\epsilon_s = o_{a.s.}(1).$$

Since $j < l$, $\epsilon_{s-(l-j)}$ is $\mathcal{F}_{s-1}$-measurable, as the $\sigma$-fields are increasing. Set $\tilde{z}_s = \kappa_j^\phi(s)\kappa_l^\phi(s)\epsilon_{s-(l-j)}$, which is $\mathcal{F}_{s-1}$-measurable by part (b) of Proposition 4.7. Then from boundedness, $\limsup \frac{1}{t}\sum_{s=1}^{t}\tilde{z}_s^2 \leq (\sup_t |\kappa^\phi(t)|)^4 (\sup_t |\epsilon_t|)^2 < \infty$, and thus from Proposition 4.2, $t^{-1}\sum_{s=1}^{t} \tilde{z}_s \epsilon_s \xrightarrow{a.s.} 0$ and therefore (4.21) holds and the proof of the Lemma is complete. □

**Lemma 4.2.** *For each $u \geq 0$, under the assumptions of Theorem 4.1,*

$$\text{(4.28)} \qquad \frac{1}{t}\sum_{s=1}^{t}\phi_s^2 = \frac{1}{t}\sum_{s=1}^{t}\left(\sum_{j=0}^{u}\kappa_j^\phi(s)\epsilon_{s-j}\right)^2 + r_1(t,u)$$

*where $\lim_u \limsup_t |r_1(t,u)| = 0$.*

*Proof.*

$$\frac{1}{t}\sum_{s=1}^{t}\phi_s^2 = \frac{1}{t}\sum_{s=1}^{t}\left(\sum_{j=0}^{\infty}\kappa_j^\phi(s)\epsilon_{s-j}\right)^2$$

$$= \frac{1}{t}\sum_{s=1}^{t}\left(\sum_{j=0}^{u}\kappa_j^\phi(s)\epsilon_{s-j}\right)^2 + \frac{2}{t}\sum_{s=1}^{t}\left(\sum_{j=0}^{u}\kappa_j^\phi(s)\epsilon_{s-j}\sum_{l=u+1}^{\infty}\kappa_l^\phi(s)\epsilon_{s-l}\right)$$

$$+ \frac{1}{t}\sum_{s=1}^{t}\left(\sum_{j=u+1}^{\infty}\kappa_j^\phi(s)\epsilon_{s-j}\right)^2.$$

Let $r_1(t,u)$ be the sum of the last two terms. Recall from Proposition 4.7 and (4.15) that $|\kappa_j^\phi(t)| \leq \kappa_j^\phi$ where $\sum_{j=0}^{\infty}\kappa_j^\phi < \infty$ a.s. From this and (D4), it follows that

$$\lim_u \limsup_t |r_1(t,u)| \leq K(\xi)\lim_u \left\{\sum_{j=0}^{u}\kappa_j^\phi \sum_{l=u+1}^{\infty}\kappa_l^\phi + \left(\sum_{j=u+1}^{\infty}\kappa_j^\phi\right)^2\right\} = 0,$$

and consequently that (4.28) holds. □

**Lemma 4.3.** *For each $u \geq 0$, under the assumptions of Theorem 4.1,*

$$\text{(4.29)} \qquad \frac{1}{t}\sum_{s=1}^{t}\left(\sum_{j=0}^{u}\kappa_j^\phi(s)\epsilon_{s-j}\right)^2 = \frac{\sigma_\epsilon^2}{t}\sum_{s=1}^{t}\sum_{j=0}^{u}\left(\kappa_j^\phi(s)\right)^2 + o_{a.s.}(1).$$

*Proof.*

$$\frac{1}{t}\sum_{s=1}^{t}\left(\sum_{j=0}^{u}\kappa_j^\phi(s)\epsilon_{s-j}\right)^2 = \frac{1}{t}\sum_{s=1}^{t}\sum_{j=0}^{u}(\kappa_j^\phi(s))^2\epsilon_{s-j}^2 + \frac{1}{t}\sum_{s=1}^{t}\sum_{j\neq l}^{u}\kappa_j^\phi(s)\epsilon_{s-j}\kappa_l^\phi(s)\epsilon_{s-l}.$$

Since $u$ is finite, by Lemma 4.1, part (c),

$$\frac{1}{t}\sum_{s=1}^{t}\sum_{j\neq l}^{u}\kappa_j^\phi(s)\epsilon_{s-j}\kappa_l^\phi(s)\epsilon_{s-l} = o_{a.s.}(1),$$

and so it remains to consider $t^{-1}\sum_{s=1}^{t}\sum_{j=0}^{u}(\kappa_j^\phi(s))^2\epsilon_{s-j}^2$. Consider the martingale difference sequence $\tilde{\epsilon}_t = \epsilon_t^2 - E[\epsilon_t^2|\mathcal{F}_{t-1}] = \epsilon_t^2 - \sigma_\epsilon^2$ (recall that $E[\epsilon_t^2|\mathcal{F}_{t-1}] = \sigma_\epsilon^2$). From (D4), $\tilde{\epsilon}_t$ is bounded a.s., hence $\sup_{-\infty<t<\infty}E[|\tilde{\epsilon}_t|^p|\mathcal{F}_{t-1}^\epsilon] < \infty$ a.s., so we can apply Proposition 4.2 to $\tilde{\epsilon}_t$.

For any $j \leq s$, consider $t^{-1}\sum_{s=1}^{t}(\kappa_j^\phi(s))^2\tilde{\epsilon}_s$. Since $\limsup t^{-1}\sum_{s=1}^{t}(\kappa_j^\phi(s))^2 < \infty$, by Proposition 4.2, $t^{-1}\sum_{s=1}^{t}(\kappa_j^\phi(s))^2\tilde{\epsilon}_s \xrightarrow{a.s.} 0$, hence, $t^{-1}\sum_{s=1}^{t}(\kappa_j^\phi(s))^2(\epsilon_s^2 - \sigma_\epsilon^2) \xrightarrow{a.s.} 0$. By an argument like that used to prove part (b) of Lemma 4.1, it follows that

$$\frac{1}{t}\sum_{s=1}^{t}(\kappa_j^\phi(s))^2\tilde{\epsilon}_{s-j} = \frac{1}{t}\sum_{s=j+1}^{t}(\kappa_j^\phi(s-j))^2\tilde{\epsilon}_{s-j} + o_{a.s.}(1),$$

$$= \frac{1}{t}\sum_{s=1}^{t-j}(\kappa_j^\phi(s))^2\tilde{\epsilon}_s + o_{a.s.}(1),$$

$$= \frac{1}{t}\sum_{s=1}^{t}(\kappa_j^\phi(s))^2\tilde{\epsilon}_s + o_{a.s.}(1) \ \text{(since } j \leq s\text{)},$$

$$= \frac{1}{t}\sum_{s=1}^{t}(\kappa_j^\phi(s))^2\tilde{\epsilon}_s + o_{a.s.}(1) = o_{a.s.}(1),$$

i.e., $t^{-1}\sum_{s=1}^{t}(\kappa_j^\phi(s))^2(\epsilon_{s-j}^2 - \sigma_\epsilon^2) \xrightarrow{a.s.} 0$ for all $j \leq u$. Finally, since $u$ is finite, $t^{-1}\sum_{s=1}^{t}\sum_{j=0}^{u}(\kappa_j^\phi(s))^2(\epsilon_{s-j}^2 - \sigma_\epsilon^2) \xrightarrow{a.s.} 0$, and (4.29) holds and the proof of the Lemma is complete. □

**Lemma 4.4.** *Under the assumptions of Theorem 4.1, for each $u \geq 0$ and $0 \leq k^* < \infty$, we have*

$$(4.30) \quad \frac{\sigma_\epsilon^2}{t}\sum_{s=1}^{t}\sum_{j=0}^{u}\left(\kappa_j^\phi(s)\right)^2 = \frac{\sigma_\epsilon^2}{t}\sum_{s=1}^{t}\sum_{j=0}^{u}\left(\sum_{l=0}^{j}(-\theta_{s+k^*})^l\sum_{p=0}^{j-l}(-\beta\theta_{s+k^*})^p\kappa_{j-l-p}\right)^2 + o_{a.s.}(1).$$

*Proof.* First suppose $k^* = 0$. Recalling from (4.12) and (4.14), for $s \geq j+1$,

$$\kappa_j^\phi(s) = \sum_{l=0}^{j}(-1)^l\kappa_{j-l}^x(s)\prod_{i=1}^{l}\theta_{s-i}$$

$$(4.31) \quad = \sum_{l=0}^{j}(-1)^l\left(\sum_{p=0}^{j-l}(-\beta)^p\kappa_{j-l-p}\prod_{r=1}^{p}\theta_{s-r}\right)\prod_{i=1}^{l}\theta_{s-i}.$$

From (3.2) and (4.31), it follows that for $s \geq j+1$,

$$\kappa_j^\phi(s) = \sum_{l=0}^{j}(-1)^l\kappa_{j-l}^x(s)\prod_{i=1}^{l}\theta_{s-i} = \sum_{l=0}^{j}(-1)^l\theta_{s-1}\kappa_{j-l}^x(s)\prod_{i=2}^{l}\theta_{s-i}$$

$$(4.32) \quad = \sum_{l=0}^{j}(-1)^l(\theta_s - s^{-1}\bar{P}_s^{-1}\phi_{s-1}e_s)\kappa_{j-l}^x(s)\prod_{i=2}^{l}\theta_{s-i}$$

$$= \sum_{l=0}^{j} (-1)^l \theta_s \kappa^x_{j-l}(s) \prod_{i=2}^{l} \theta_{s-i}$$

$$- \sum_{l=0}^{j} (-1)^l s^{-1} \bar{P}_s^{-1} \phi_{s-1} e_s \kappa^x_{j-l}(s) \prod_{i=2}^{l} \theta_{s-i}.$$

Next, taking the square of (4.32), we obtain

(4.33)
$$(\kappa^\phi_j(s))^2 = \left( \sum_{l=0}^{j} (-1)^l \theta_s \kappa^x_{j-l}(s) \prod_{i=2}^{l} \theta_{s-i} \right)^2$$
$$- 2 \sum_{l=0}^{j} (-1)^l$$
$$\times \sum_{m=0}^{j} \left\{ (-1)^m (\theta_s \kappa^x_{j-l}(s) \left( \prod_{i=2}^{l} \theta_{s-i} \right) s^{-1} \bar{P}_s^{-1} \phi_{s-1} e_s) \kappa^x_{j-m}(s) \prod_{p=2}^{m} \theta_{s-p} \right\}$$
$$+ \left( \sum_{l=0}^{j} (-1)^l s^{-1} \bar{P}_s^{-1} \phi_{s-1} e_s \kappa^x_{j-l}(s) \prod_{i=2}^{l} \theta_{s-i} \right)^2,$$

and from the boundedness of $\theta_t$ and an argument like that used to prove (4.17), it follows that

(4.34)
$$\frac{\sigma_\epsilon^2}{t} \sum_{s=j+1}^{t} (\kappa^\phi_j(s))^2 = \frac{\sigma_\epsilon^2}{t} \sum_{s=j+1}^{t} \left( \sum_{l=0}^{j} (-1)^l \theta_s \kappa^x_{j-l}(s) \prod_{i=2}^{l} \theta_{s-i} \right)^2 + o_{a.s.}(1).$$

Consider next the r.h.s. of (4.34). From (3.2),

$$\sum_{l=0}^{j} (-1)^l \theta_s \kappa^x_{j-l}(s) \prod_{i=2}^{l} \theta_{s-i}$$
$$= \sum_{l=0}^{j} (-1)^l \theta_s \theta_{s-2} \kappa^x_{j-l}(s) \prod_{i=3}^{l} \theta_{s-i}$$
$$= \sum_{l=0}^{j} (-1)^l \theta_s \left( \theta_s - \sum_{m=0}^{1} (s-m)^{-1} \bar{P}_{s-m}^{-1} \phi_{s-m-1} e_{s-m} \right) \kappa^x_{j-l}(s) \prod_{i=3}^{l} \theta_{s-i}$$
$$= \sum_{l=0}^{j} (-1)^l \theta_s^2 \kappa^x_{j-l}(s) \prod_{i=3}^{l} \theta_{s-i}$$
$$- \sum_{l=0}^{j} (-1)^l \theta_s \sum_{m=0}^{1} (s-m)^{-1} \bar{P}_{s-m}^{-1} \phi_{s-m-1} e_{s-m} \kappa^x_{j-l}(s) \prod_{i=3}^{l} \theta_{s-i}.$$

Therefore, again from boundedness of $\theta_t$ and an argument like that used to prove (4.17),

$$\frac{\sigma_\epsilon^2}{t} \sum_{s=j+1}^{t} (\kappa^\phi_j(s))^2 = \frac{\sigma_\epsilon^2}{t} \sum_{s=j+1}^{t} \left( \sum_{l=0}^{j} (-1)^l \theta_s^2 \kappa^x_{j-l}(s) \prod_{i=3}^{l} \theta_{s-i} \right)^2 + o_{a.s.}(1).$$

Applying the argument $l - 2$ additional times, it follows that

$$(4.35) \quad \frac{\sigma_\epsilon^2}{t} \sum_{s=j+1}^{t} (\kappa_j^\phi(s))^2 = \frac{\sigma_\epsilon^2}{t} \sum_{s=j+1}^{t} \left( \sum_{l=0}^{j} (-\theta_s)^l \kappa_{j-l}^x(s) \right)^2 + o_{a.s.}(1).$$

Next working on the r.h.s. of (4.35), from (4.12):

$$\sum_{l=0}^{j} (-\theta_s)^l \kappa_{j-l}^x(s) = \sum_{l=0}^{j} (-\theta_s)^l \left( \sum_{p=0}^{j-l} (-\beta)^p \kappa_{j-l-p} \prod_{r=1}^{p} \theta_{s-r} \right)$$

$$= \sum_{l=0}^{j} (-\theta_s)^l \left( \sum_{p=0}^{j-l} (-\beta)^p (\theta_s - s^{-1} \bar{P}_s^{-1} \phi_{s-1} e_s) \kappa_{j-l-p} \prod_{r=2}^{p} \theta_{s-r} \right)$$

$$= \sum_{l=0}^{j} (-\theta_s)^l \left( \sum_{p=0}^{j-l} (-\beta)^p \theta_s \kappa_{j-l-p} \prod_{r=2}^{p} \theta_{s-r} \right)$$

$$- \sum_{l=0}^{j} (-\theta_s)^l \left( \sum_{p=0}^{j-l} (-\beta)^p s^{-1} \bar{P}_s^{-1} \phi_{s-1} e_s \kappa_{j-l-p} \prod_{r=2}^{p} \theta_{s-r} \right).$$

Hence,

$$\frac{\sigma_\epsilon^2}{t} \sum_{s=j+1}^{t} \left( \sum_{l=0}^{j} (-\theta_s)^l \kappa_{j-l}^x(s) \right)^2$$

$$= \frac{\sigma_\epsilon^2}{t} \sum_{s=j+1}^{t} \left( \sum_{l=0}^{j} (-\theta_s)^l \sum_{p=0}^{j-l} (-\beta)^p \theta_s \kappa_{j-l-p} \prod_{r=2}^{p} \theta_{s-r} \right)^2 + o_{a.s.}(1).$$

Applying the argument $p - 1$ additional times, it follows that

$$\frac{\sigma_\epsilon^2}{t} \sum_{s=j+1}^{t} \left( \sum_{l=0}^{j} (-\theta_s)^l \kappa_{j-l}^x(s) \right)^2 = \frac{\sigma_\epsilon^2}{t} \sum_{s=j+1}^{t} \left( \sum_{l=0}^{j} (-\theta_s)^l \sum_{p=0}^{j-l} (-\beta \theta_s)^p \kappa_{j-l-p} \right)^2$$
$$+ o_{a.s.}(1),$$

and, since $j$ is finite,

$$\frac{\sigma_\epsilon^2}{t} \sum_{s=1}^{t} (\kappa_j^\phi(s))^2 = \frac{\sigma_\epsilon^2}{t} \sum_{s=j+1}^{t} \left( \sum_{l=0}^{j} (-\theta_s)^l \sum_{p=0}^{j-l} (-\beta \theta_s)^p \kappa_{j-l-p} \right)^2 + o_{a.s.}(1).$$

Finally, since $u$ is finite, (4.30) follows for $k^* = 0$.

From (3.2), for any finite $k^* > 0$,

$$\theta_{s+k^*} = \theta_s + \sum_{r=0}^{k^*-1} (s + k^* - r)^{-1} \bar{P}_{s+k^*-r}^{-1} \phi_{s+k^*-r-1} e_{s+k^*-r}.$$

Set

$$\lambda(s, k^*) = \sum_{r=0}^{k^*-1} (s + k^* - r)^{-1} \bar{P}_{s+k^*-r}^{-1} \phi_{s+k^*-r-1} e_{s+k^*-r}.$$

For every integer $l \geq 0$, the binomial formula yields

$$\theta_{s+k^*}^l = \theta_s^l + \binom{l}{1}\theta_{s+k^*}^{l-1}\lambda(s,k^*) + \cdots + \binom{l}{l-1}\theta_{s+k^*}\lambda^{l-1}(s,k^*) + \lambda^l(s,k^*).$$

Substituting this result into the r.h.s. of

$$\frac{\sigma_\epsilon^2}{t}\sum_{s=1}^t (\kappa_j^\phi(s))^2 = \frac{\sigma_\epsilon^2}{t}\sum_{s=j+1}^t \left(\sum_{l=0}^j (-\theta_{s+k^*})^l \sum_{p=0}^{j-l}(-\beta\theta_{s+k^*})^p \kappa_{j-l-p}\right)^2 + o_{a.s.}(1),$$

which follows from Lemma 4.1, and noting that each resulting term involving $\lambda(s,k^*)$ is $o_{a.s.}(1)$ by (4.18), the proof of (4.30) and of the Lemma is reduced to the result just established for $k^* = 0$. $\square$

**Lemma 4.5.** *Under the assumptions of Theorem 4.1, for any finite $u$,*

(4.36) $$\frac{\sigma_\epsilon^2}{t}\sum_{s=2}^t \sum_{j=1}^u \kappa_j^e(s)\kappa_{j-1}^\phi(s-1)\epsilon_{s-j}^2 = \frac{\sigma_\epsilon^2}{t}\sum_{s=2}^t \sum_{j=1}^u \kappa_j^e(s)\kappa_{j-1}^\phi(s-1) + o_{a.s.}(1).$$

*Proof.* Since $u$ is finite and for any finite $j \leq u$, $\limsup t^{-1}\sigma_\epsilon^2 |\sum_{s=2}^t \kappa_j^e(s) \times \kappa_{j-1}^\phi(s-1)| < \infty$, then the result follows by an argument similar to that used to prove part (c) of Lemma 4.1. $\square$

**Proposition 4.8.** *Under the assumptions of Theorem 4.1, the sequence $\{\hat{\theta}_t\}$ defined by (4.1) satisfies $\theta_t - \hat{\theta}_t = o_{a.s.}(1)$.*

*Proof.* For simplicity, first assume that $k^* = 0$; i.e., $|\theta_t| \leq K^* < 1$ for all $t$. From the results of Proposition 4.7 and Lemmas 4.2–4.4, for any $u < \infty$:

$$\frac{1}{t}\sum_{s=1}^t \phi_s^2 = \frac{1}{t}\sum_{s=1}^t \left(\sum_{j=0}^u \kappa_j^\phi(s)\epsilon_{s-j}\right)^2 + r_1(t,u) \quad \text{(Lemma 4.2)}$$

$$= \frac{\sigma_\epsilon^2}{t}\sum_{s=1}^t \sum_{j=0}^u \left(\kappa_j^\phi(s)\right)^2 + o_{a.s.}(1) + r_1(t,u) \quad \text{(Lemma 4.3)}$$

$$= \frac{\sigma_\epsilon^2}{t}\sum_{s=1}^t \sum_{j=0}^u \left(\sum_{l=0}^j (-1)^l \prod_{i=1}^l \theta_{s-i} \sum_{p=0}^{j-l}(-\beta)^p \kappa_{j-l-p} \prod_{r=1}^p \theta_{s-r}\right)^2 + o_{a.s.}(1)$$
$$+ r_1(t,u)$$

$$= \frac{\sigma_\epsilon^2}{t}\sum_{s=1}^t \sum_{j=0}^u \left(\sum_{l=0}^j (-\theta_s)^l \sum_{p=0}^{j-l}(-\beta\theta_s)^p \kappa_{j-l-p}\right)^2 + o_{a.s.}(1)$$
$$+ r_1(t,u) \quad \text{(Lemma 4.4)}$$

where $\lim_u \limsup_t |r_1(t,u)| = 0$. By (2.3), Parseval's relation and convolution [22, pp. 61-66], it follows that

$$\frac{\sigma_\epsilon^2}{t}\sum_{s=1}^t \sum_{j=0}^\infty \left(\sum_{l=0}^j (-\theta_s)^l \sum_{p=0}^{j-l}(-\beta\theta_s)^p \kappa_{j-l-p}\right)^2$$

$$= \frac{1}{t}\sum_{s=1}^t \int_{-\pi}^\pi \frac{1}{|(1+\theta_s e^{i\omega})(1+\beta\theta_s e^{i\omega})|^2} g_y(\omega)d\omega,$$

and so

$$\frac{1}{t}\sum_{s=1}^{t}\phi_s^2 = \frac{1}{t}\sum_{s=1}^{t}\int_{-\pi}^{\pi}\frac{1}{|(1+\theta_s e^{i\omega})(1+\beta\theta_s e^{i\omega})|^2}g_y(\omega)d\omega + r_2(t,u)$$

where

$$r_2(t,u) = r_1(t,u) + \frac{\sigma_\epsilon^2}{t}\sum_{s=1}^{t}\sum_{j=u+1}^{\infty}\left(\sum_{l=0}^{j}(-\theta_s)^l\sum_{p=0}^{j-l}(-\beta\theta_s)^p\kappa_{j-l-p}\right)^2.$$

Since $|\theta_t| \leq K^* < 1$, it follows that

$$\frac{\sigma_\epsilon^2}{t}\sum_{s=1}^{t}\sum_{j=u+1}^{\infty}\left(\sum_{l=0}^{j}(-\theta_s)^l\sum_{p=0}^{j-l}(-\beta\theta_s)^p\kappa_{j-l-p}\right)^2$$

$$\leq \frac{\sigma_\epsilon^2}{t}\sum_{j=u+1}^{\infty}\left(\sum_{l=0}^{j}(K^*)^l\sum_{p=0}^{j-l}(\beta K^*)^p\kappa_{j-l-p}\right)^2 \overset{u\to\infty}{\longrightarrow} 0$$

because $\sum_{j=0}^{\infty}\left(\sum_{l=0}^{j}(K^*)^l\sum_{p=0}^{j-l}(\beta K^*)^p\kappa_{j-l-p}\right)^2 < \infty$. Hence,

$$\lim_{u}\limsup_{t}\frac{\sigma_\epsilon^2}{t}\sum_{s=1}^{t}\sum_{j=u+1}^{\infty}\left(\sum_{l=0}^{j}(-\theta_s)^l\sum_{p=0}^{j-l}(-\beta\theta_s)^p\kappa_{j-l-p}\right)^2 = 0,$$

and consequently, $\lim_u \limsup_t |r_2(t,u)| = 0$. It follows that

$$(4.37) \qquad \left|\frac{1}{t}\sum_{s=1}^{t}\phi_s^2 - \frac{1}{t}\sum_{s=1}^{t}\int_{-\pi}^{\pi}\frac{1}{|(1+\theta_s e^{i\omega})(1+\beta\theta_s e^{i\omega})|^2}g_y(\omega)d\omega\right| \overset{a.s.}{\longrightarrow} 0.$$

Next, for $s \geq 2$ we use $\kappa_j^z(s) = \kappa_j^e(s) + \theta_{s-1}\kappa_{j-1}^\phi(s-1)$ from (4.16) to obtain that for any $u < \infty$,

$$\frac{1}{t}\sum_{s=2}^{t}z_s\phi_{s-1} = \frac{1}{t}\sum_{s=2}^{t}\sum_{j=0}^{\infty}\kappa_j^z(s)\epsilon_{s-j}\sum_{l=0}^{\infty}\kappa_l^\phi(s-1)\epsilon_{s-1-l}$$

$$= \frac{1}{t}\sum_{s=2}^{t}\sum_{j=0}^{u}\kappa_j^z(s)\epsilon_{s-j}\sum_{l=0}^{u}\kappa_l^\phi(s-1)\epsilon_{s-1-l} + r_3(t,u)$$

$$= \frac{1}{t}\sum_{s=2}^{t}\sum_{j=1}^{u}\kappa_j^z(s)\kappa_{j-1}^\phi(s-1)\epsilon_{s-j}^2 + o_{a.s.}(1) + r_3(t,u)$$

$$= \frac{1}{t}\sum_{s=2}^{t}\sum_{j=1}^{u}\left\{\left(\kappa_j^e(s) + \theta_{s-1}\kappa_{j-1}^\phi(s-1)\right)\kappa_{j-1}^\phi(s-1)\right\}\epsilon_{s-j}^2$$

$$\quad + o_{a.s.}(1) + r_3(t,u)$$

$$= \frac{1}{t}\sum_{s=2}^{t}\sum_{j=1}^{u}\kappa_j^e(s)\kappa_{j-1}^\phi(s-1)\epsilon_{s-j}^2$$

$$(4.38) \qquad + \frac{1}{t}\sum_{s=2}^{t}\sum_{j=1}^{u}\theta_{s-1}\left(\kappa_{j-1}^\phi(s-1)\right)^2\epsilon_{s-j}^2 + o_{a.s.}(1) + r_3(t,u)$$

$$= \frac{\sigma_\epsilon^2}{t} \sum_{s=2}^{t} \sum_{j=1}^{u} \kappa_j^e(s) \kappa_{j-1}^\phi(s-1) \text{ (Lemma 4.5)}$$

$$+ \frac{\sigma_\epsilon^2}{t} \sum_{s=2}^{t} \sum_{j=1}^{u} \theta_{s-1} \left(\kappa_{j-1}^\phi(s-1)\right)^2 + o_{a.s.}(1) + r_3(t,u),$$

where

$$r_3(t,u) = \frac{1}{t} \sum_{s=2}^{t} \left( \sum_{j=u+1}^{\infty} \kappa_j^z(s)\epsilon_{s-j} \sum_{l=0}^{u} \kappa_l^\phi(s-1)\epsilon_{s-1-l} \right)$$

$$+ \frac{1}{t} \sum_{s=2}^{t} \left( \sum_{j=0}^{u} \kappa_j^z(s)\epsilon_{s-j} \sum_{l=u+1}^{\infty} \kappa_l^\phi(s-1)\epsilon_{s-1-l} \right)$$

$$+ \frac{1}{t} \sum_{s=2}^{t} \left( \sum_{j=u+1}^{\infty} \kappa_j^z(s)\epsilon_{s-j} \sum_{l=u+1}^{\infty} \kappa_l^\phi(s-1)\epsilon_{s-1-l} \right).$$

By an argument similar to that applied to $r_1(t,u)$ in the proof of Lemma 4.2, one obtains $\lim_u \limsup_t |r_3(t,u)| = 0$. As shown above, the second term of (4.38), $t^{-1}\sigma_\epsilon^2 \sum_{s=2}^{t} \sum_{j=1}^{u} \theta_{s-1} \left(\kappa_{j-1}^\phi(s-1)\right)^2$, is equal to

$$\frac{1}{t} \sum_{s=1}^{t} \int_{-\pi}^{\pi} \frac{\theta_s}{|(1+\theta_s e^{i\omega})(1+\beta\theta_s e^{i\omega})|^2} g_y(\omega) d\omega + r_4(t,u)$$

with

$$r_4(t,u) = t^{-1}\sigma_\epsilon^2 \sum_{s=2}^{t} \sum_{j=u+1}^{\infty} \theta_{s-1} \left(\kappa_{j-1}^\phi(s-1)\right)^2 + o_{a.s.}(1)$$

and $\lim_u \limsup_t |r_4(t,u)| = 0$. Hence, it remains to consider the first term of (4.38). From (4.11) and (4.31),

$$\frac{\sigma_\epsilon^2}{t} \sum_{s=2}^{t} \sum_{j=1}^{u} \kappa_j^e(s) \kappa_{j-1}^\phi(s-1)$$

$$= \frac{\sigma_\epsilon^2}{t} \sum_{s=2}^{t} \sum_{j=1}^{u} \left\{ \sum_{l=0}^{j} (-1)^l \kappa_{j-l} \prod_{i=1}^{l} \theta_{s-i} \right\}$$

$$\times \left\{ \sum_{m=0}^{j-1} (-1)^m \prod_{p=1}^{m} \theta_{s-1-p} \left( \sum_{n=0}^{j-l-1} (-\beta)^n \kappa_{j-l-1-n} \prod_{r=1}^{n} \theta_{s-1-r} \right) \right\}$$

$$= \frac{\sigma_\epsilon^2}{t} \sum_{s=2}^{t} \sum_{j=1}^{u} \left\{ \left( \sum_{l=0}^{j} (-\theta_s)^l \kappa_{j-l} \right) \left( \sum_{m=0}^{j-1} (-\theta_s)^m \sum_{n=0}^{j-l-1} (-\beta\theta_s)^n \kappa_{j-l-1-n} \right) \right\}$$

$$+ o_{a.s.}(1),$$

and, since again by (2.3), Parseval's relation and convolution,

$$\sigma_\epsilon^2 \sum_{j=1}^{\infty} \left\{ \left( \sum_{l=0}^{j} (-\theta_s)^l \kappa_{j-l} \right) \left( \sum_{m=0}^{j-1} (-\theta_s)^m \sum_{n=0}^{j-l-1} (-\beta\theta_s)^n \kappa_{j-l-1-n} \right) \right\}$$

$$= \int_{-\pi}^{\pi} \frac{1}{(1+\theta_s e^{-i\omega})} \frac{e^{i\omega}}{(1+\theta_s e^{i\omega})(1+\beta\theta_s e^{i\omega})} g_y(\omega) d\omega,$$

the first term of (4.38), $t^{-1}\sigma_\epsilon^2 \sum_{s=2}^t \sum_{j=1}^u \kappa_j^e(s)\kappa_{j-1}^\phi(s-1)$, is equal to

$$\frac{1}{t}\sum_{s=1}^t \frac{1}{(1+\theta_s e^{-i\omega})}\frac{e^{i\omega}}{(1+\theta_s e^{i\omega})(1+\beta\theta_s e^{i\omega})} g_y(\omega)d\omega + r_5(t,u)$$

with

$$r_5(t,u) = \frac{\sigma_\epsilon^2}{t}\sum_{s=1}^t \sum_{j=u+1}^\infty \left\{\left(\sum_{l=0}^j (-\theta_s)^l \kappa_{j-l}\right) \right.$$
$$\left. \times \left(\sum_{m=0}^{j-1}(-\theta_s)^m \sum_{n=0}^{j-l-1}(-\beta\theta_s)^n \kappa_{j-l-1-n}\right)\right\} + o_{a.s.}(1).$$

An argument like that applied to $r_2(t,u)$ yields $\lim_u \limsup_t |r_5(t,u)| = 0$. Further, since

$$\int_{-\pi}^{\pi} \frac{1}{(1+\theta_s e^{-i\omega})}\frac{e^{i\omega}}{(1+\theta_s e^{i\omega})(1+\beta\theta_s e^{i\omega})} g_y(\omega)d\omega$$
$$+ \int_{-\pi}^{\pi} \frac{\theta_s}{|(1+\theta_s e^{i\omega})(1+\beta\theta_s e^{i\omega})|^2} g_y(\omega)d\omega$$
$$= \int_{-\pi}^{\pi} \frac{e^{i\omega}(1+\beta\theta_s e^{-i\omega}) + \theta_s}{|(1+\theta_s e^{i\omega})(1+\beta\theta_s e^{i\omega})|^2} g_y(\omega)d\omega$$
$$= \int_{-\pi}^{\pi} \frac{e^{i\omega}+(1+\beta)\theta_s}{|(1+\theta_s e^{i\omega})(1+\beta\theta_s e^{i\omega})|^2} g_y(\omega)d\omega,$$

we obtain

$$(4.39) \quad \left|\frac{1}{t}\sum_{s=2}^t z_s \phi_{s-1} - \frac{1}{t}\sum_{s=1}^t \int_{-\pi}^{\pi} \frac{e^{i\omega}+(1+\beta)\theta_s}{|(1+\theta_s e^{i\omega})(1+\beta\theta_s e^{i\omega})|^2} g_y(\omega)d\omega\right| \xrightarrow{a.s.} 0.$$

Combining (4.37) and (4.39), it follows from (3.5) and (3.16) that $\theta_t - \hat{\theta}_t = o_{a.s.}(1)$ where $\hat{\theta}_t$ is given by (4.1). Finally, since $k^*$ is finite, an argument similar to that used for Lemma 4.4 can be applied to show that (4.1) holds for the general case $k^* > 0$, completing the proofs of the proposition and part (a) of Theorem 4.1. □

**Proposition 4.9.** *Under the assumptions of Theorem 4.1 and with $t_0 = t_0(\xi)$ as in (a) of the Theorem, $\hat{\theta}_t$ defined by (4.1) satisfies the conditions of Proposition 4.6 for $\tilde{\Theta} = \Theta = (-1,1)$ for a Robbins-Monro recursion with $f(\theta) = f(\theta, \beta)$ as in (1.4).*

*Proof.* For $t \geq 2$, set

$$\tilde{P}_t = \sum_{s=1}^t \int_{-\pi}^{\pi} \frac{1}{|(1+\theta_{s+k^*} e^{i\omega})(1+\beta\theta_{s+k^*} e^{i\omega})|^2} g_y(\omega)d\omega.$$

Then from (4.1),

$$\hat{\theta}_t = \tilde{P}_t^{-1} \sum_{s=1}^t \int_{-\pi}^{\pi} \frac{e^{i\omega}+(1+\beta)\theta_{s+k^*}}{|(1+\theta_{s+k^*} e^{i\omega})(1+\beta\theta_{s+k^*} e^{i\omega})|^2} g_y(\omega)d\omega$$
$$= \tilde{P}_t^{-1} \left\{\sum_{s=1}^{t-1} \int_{-\pi}^{\pi} \frac{e^{i\omega}+(1+\beta)\theta_{s+k^*}}{|(1+\theta_{s+k^*} e^{i\omega})(1+\beta\theta_{s+k^*} e^{i\omega})|^2} g_y(\omega)d\omega \right.$$
$$(4.40) \qquad \left. + \int_{-\pi}^{\pi} \frac{e^{i\omega}+(1+\beta)\theta_{t+k^*}}{|(1+\theta_{t+k^*} e^{i\omega})(1+\beta\theta_{t+k^*} e^{i\omega})|^2} g_y(\omega)d\omega\right\}$$

$$= \tilde{P}_t^{-1} \left\{ \tilde{P}_{t-1} \hat{\theta}_{t-1} + \int_{-\pi}^{\pi} \frac{e^{i\omega} + (1+\beta)\theta_{t+k^*}}{|(1+\theta_{t+k^*}e^{i\omega})(1+\beta\theta_{t+k^*}e^{i\omega})|^2} g_y(\omega) d\omega \right\}$$

$$= \tilde{P}_t^{-1} \left\{ \left( \tilde{P}_t - \int_{-\pi}^{\pi} \frac{1}{|(1+\theta_{t+k^*}e^{i\omega})(1+\beta\theta_{t+k^*}e^{i\omega})|^2} g_y(\omega) d\omega \right) \hat{\theta}_{t-1} \right.$$
$$\left. + \int_{-\pi}^{\pi} \frac{e^{i\omega} + (1+\beta)\theta_{t+k^*}}{|(1+\theta_{t+k^*}e^{i\omega})(1+\beta\theta_{t+k^*}e^{i\omega})|^2} g_y(\omega) d\omega \right\}$$

$$= \hat{\theta}_{t-1} + \tilde{P}_t^{-1} \left\{ \int_{-\pi}^{\pi} \frac{e^{i\omega} + (1+\beta)\theta_{t+k^*}}{|(1+\theta_{t+k^*}e^{i\omega})(1+\beta\theta_{t+k^*}e^{i\omega})|^2} g_y(\omega) d\omega \right.$$
$$\left. - \int_{-\pi}^{\pi} \frac{\hat{\theta}_{t-1}}{|(1+\theta_{t+k^*}e^{i\omega})(1+\beta\theta_{t+k^*}e^{i\omega})|^2} g_y(\omega) d\omega \right\}$$

$$= \hat{\theta}_{t-1} + \tilde{P}_t^{-1} \int_{-\pi}^{\pi} \frac{e^{i\omega} + \beta\theta_{t+k^*}}{|(1+\theta_{t+k^*}e^{i\omega})(1+\beta\theta_{t+k^*}e^{i\omega})|^2} g_y(\omega) d\omega$$
$$+ \tilde{P}_t^{-1}(\theta_{t+k^*} - \hat{\theta}_{t-1}) \int_{-\pi}^{\pi} \frac{1}{|(1+\theta_{t+k^*}e^{i\omega})(1+\beta\theta_{t+k^*}e^{i\omega})|^2} g_y(\omega) d\omega$$

$$= \hat{\theta}_{t-1} - \tilde{P}_t^{-1} f(\theta_{t+k^*}, \beta)$$
$$+ \tilde{P}_t^{-1}(\theta_{t+k^*} - \hat{\theta}_{t-1}) \int_{-\pi}^{\pi} \frac{1}{|(1+\theta_{t+k^*}e^{i\omega})(1+\beta\theta_{t+k^*}e^{i\omega})|^2} g_y(\omega) d\omega$$

$$= \hat{\theta}_{t-1} - \tilde{P}_t^{-1} f(\hat{\theta}_{t-1}, \beta) + \tilde{P}_t^{-1} \left( f(\hat{\theta}_{t-1}, \beta) - f(\theta_{t+k^*}, \beta) \right)$$
$$+ \tilde{P}_t^{-1}(\theta_{t+k^*} - \hat{\theta}_{t-1}) \int_{-\pi}^{\pi} \frac{1}{|(1+\theta_{t+k^*}e^{i\omega})(1+\beta\theta_{t+k^*}e^{i\omega})|^2} g_y(\omega) d\omega$$

$$= \hat{\theta}_{t-1} - \delta_t f(\hat{\theta}_{t-1}, \beta) + \delta_t \gamma_t,$$

where, for $t \geq 2$,

$$(4.41) \quad \delta_t = \tilde{P}_t^{-1} = \frac{1}{t} \left[ \frac{1}{t} \sum_{s=1}^{t} \int_{-\pi}^{\pi} \frac{1}{|(1+\theta_{s+k^*}e^{i\omega})(1+\beta\theta_{s+k^*}e^{i\omega})|^2} g_y(\omega) d\omega \right]^{-1},$$

and, for $t \geq t_0 + 1$ with $t_0$ as in (a) of Theorem 4.1 (which guarantees that $f(\hat{\theta}_{t-1}, \beta)$ below is finite),

$$(4.42) \quad \begin{aligned} \gamma_t &= \left( f(\hat{\theta}_{t-1}, \beta) - f(\theta_{t+k^*}, \beta) \right) \\ &+ (\theta_{t+k^*} - \hat{\theta}_{t-1}) \int_{-\pi}^{\pi} \frac{1}{|(1+\theta_{t+k^*}e^{i\omega})(1+\beta\theta_{t+k^*}e^{i\omega})|^2} g_y(\omega) d\omega. \end{aligned}$$

For $|\theta| \leq K^* < 1$, it follows from (3.16) and (4.41), there exist finite, positive $\tilde{K}_1(\xi) \leq \tilde{K}_2(\xi)$ such that $0 < \tilde{K}_1(\xi) \leq t\delta_t \leq \tilde{K}_2(\xi) < \infty$. From this, it follows that $\delta_t \xrightarrow{a.s.} 0$ and $\sum_{s=1}^{t} \delta_s \geq \tilde{K}_1 \sum_{s=1}^{t} k^{-1} \to \infty$. Next since $k^*$ is finite, it follows from $\theta_{t-1} - \hat{\theta}_{t-1} = o_{a.s.}(1)$ (Proposition 4.7) that $\theta_{t+k^*} - \hat{\theta}_{t-1} = o_{a.s.}(1)$ and $f(\theta_{t+k^*}, \beta) - f(\hat{\theta}_{t-1}, \beta) = o_{a.s.}(1)$. Hence, $\gamma_t \xrightarrow{a.s.} 0$. The definition of $t_0$ in (a) of Theorem 4.1, guarantees that the remaining condition of Proposition 4.6 is satisfied, so the proposition is proved. □

We can now complete the proof of Theorem 4.1. By Proposition 4.6, $\hat{\theta}_t \xrightarrow{a.s.} \Theta_0^\beta$ and therefore also $\theta_t \xrightarrow{a.s.} \Theta_0^\beta$, which is compact by Proposition 4.4. Further, if $y_t$ is an invertible ARMA process, then by Proposition 4.5, the set $\Theta_0^\beta$ is finite and Proposition 4.6 shows $\theta_t$ converges on almost every realization to one of the finitely many $\theta \in \Theta_0^\beta$. Consequently, on the probability one event on which $\theta_t$ converges, its limit is a random variable $\theta$ with finitely many values. On the complementary event, $\theta$ can be defined to have any fixed value. This completes the proof of part (b) and with it the proof of the Theorem.

## 5. Discussion

The results obtained here provide a rigorous foundation for analyzing PLR and $RML_2$ for MA(1) models. An important conclusion from our results is that under misspecification, generally only $RML_2$ (i.e., the general algorithm with $\beta = 1$), not the simpler and more frequently considered PLR algorithm, can produce optimal coefficient estimates in the limit. In [5], Theorem 4.1 is applied to address convergence of PLR and convergence of $RML_2$ with a specific monitoring and modification scheme to ensure that iterates satisfy $|\theta_t| \leq K^* < 1$. In [5] we also provide a set of examples that show that the limits of $\theta_t$ from PLR and $RML_2$ can differ.

Ideas and techniques from the analysis of Hannan [12] of $RML_2$ for ARMA models played a key role in our analysis, particularly the idea of approximating the recursive algorithm's sequence by a sequence made more analyzable, replacing certain terms by their expected values, and replacing terms in an expression by finitely lagged values, as in our (4.20), so that martingale results like Propositions 4.1 and 4.2 can be applied. However we note that, because of a neglected $o_{a.s.}(1)$ term that depends on $\theta_t$, Hannan did not actually establish that his auxiliary sequence, which we denote by $\tilde{\theta}_t$ to distinguish it from our $\hat{\theta}_t$, satisfies his (nonstandard) recursion scheme. Also, the convergence analysis he indicates for $\tilde{\theta}_t$, if its details could be verified, would only establish that the *limit inferior* of $\min_{\theta \in \Theta_0^1} |\tilde{\theta}_t - \theta|$ is zero a.s., see p. 773 of Hannan [12]. The stronger result with the *limit* is needed to establish convergence of the original recursive sequence to $\Theta_0^1$. More information about problems we encountered with analyses in Hannan [12] can be found in [4, Appendix E].

The approximating sequence technique is similar to the Ordinary Differential Equation (ODE) method independently developed by Ljung [20] and Kushner [14, 15]. Specifically, the ODE method is a technique for providing asymptotic analysis of a time series (discrete stochastic process) via a deterministic continuous time stability analysis of a set of ODEs. For example, from the ODE method, Ljung makes convergence assertions for both PLR and $RML_2$ for ARMAX models including in the incorrect model situation [21]. Like Hannan, however, the analysis is incomplete. In the rigorous treatment of the ODE method presented by Benveniste et al [1] only the correct model situation is considered. Their results, however, do not apply to PLR or $RML_2$ [4, pp.65-67].

Clearly, the boundedness assumption (D4) is restrictive but it is typical in convergence analyses like ours. For example, boundedness is an explicit assumption in the deep correct model results obtained by Lai and Ying [17, equation (1.3)] as well as Ljung's ODE method assertions [21, condition S2, p.191] and is also required in the treatment by Benveniste et al. in which $\theta_t$ is assumed to be bounded to obtain verifiable conditions to prove asymptotic results [1, Theorem 15 and Corollary 16, p.238].

Finally, it is likely that Theorem 4.1 is generalizable to higher order moving average models and quite possibly ARMA models. However, to obtain convergence results, a multidimensional parameter vector $\theta$ version of Proposition 4.6 is needed.

The proof of Theorem 2.2.2 of [7] seems to provide the needed result if it can be shown that an appropriate Liapounov function exists for the vector-valued $f(\theta, \beta)$ associated with multidimensional $\theta$ for $0 \leq \beta < 1$. The generalization of the $\bar{L}(\theta)$ with vector $\theta$ provides the Liapounov function for the case $\beta = 1$.

## Acknowledgment

The authors gratefully acknowledge the detailed review and insightful comments and suggestions of the referee.